\documentclass[lettersize,journal,twoside]{IEEEtran}
\usepackage{amsmath}
\usepackage{amsfonts}
\usepackage{etoolbox}
\patchcmd{\subequations}{}%
{}{}{}

\usepackage[ruled, algonl, vlined]{algorithm2e}
\SetAlFnt{\small}
\SetAlCapFnt{\small}
\SetAlCapNameFnt{\small}

\usepackage{array}
\usepackage{booktabs}
\usepackage[caption=false,font=normalsize,labelfont=sf,textfont=sf]{subfig}
\usepackage{cite}
\usepackage{anyfontsize}

\usepackage{enumitem}
\usepackage{float}
\usepackage{graphicx}
\usepackage{multirow}
\usepackage{nomencl}
\usepackage{stfloats}
\usepackage{setspace}
\usepackage{textcomp}
\usepackage{tabularx}
\usepackage{url}
\usepackage{verbatim}
\usepackage{makecell}
\usepackage{dcolumn}
\usepackage{pifont}
\usepackage{setspace}
\usepackage{color}
\usepackage{hyperref}  
\hypersetup{colorlinks = true,
            linkcolor  = red,
            citecolor  = red,
            urlcolor   = red}
\definecolor{CBlue}{RGB}{20, 80, 163}

\newtheorem{proposition}{Proposition}

\begin{document}

\allowdisplaybreaks

\title{Towards Improving Unit Commitment Economics:\\ An Add-On Tailor for\\ Renewable Energy and Reserve Predictions}

\author{Xianbang Chen,~\IEEEmembership{Student Member,~IEEE,}~Yikui Liu,~\IEEEmembership{Member,~IEEE,}~Lei Wu,~\IEEEmembership{Fellow,~IEEE}

\thanks{

Manuscript received 4 October 2023; revised 4 February 2024 and 12 May 2024; accepted 28 June 2024. Date of publication DD MMMM YYYY; date of current version DD MMMM YYYY. This work was supported in part by the Army Combat Capabilities Development Command (CCDC) project on Resiliency of Energy Resources and Supply Chain for the Industrial Base and the PSEG Foundation gift. Paper no. TSTE-xxxxx-xxxx. \textit{(Corresponding author: Yikui Liu.)}

X. Chen and L. Wu are with the ECE Department, Stevens Institute of Technology, Hoboken, NJ, 07030 USA. (e-mail: xchen130@stevens.edu and lei.wu@stevens.edu).

Y. Liu is with the Electrical Engineering College, Sichuan University, Chengdu, 610017 China. (e-mail: yikuiliu89@outlook.com).

Color versions of one or more of the figures in this paper are available online at http://ieeexplore.ieee.org.

Digital Object Identifier xx.xxxx/TSTE.xxxx.xxxxxxx
} }

\markboth{IEEE TRANSACTIONS ON SUSTAINABLE ENERGY, Accepted}{CHEN \MakeLowercase{\textit{et al.}}: Add-On Tailor for Renewable Energy and Reserve Predictions}


\maketitle

\begin{abstract}
Generally, day-ahead unit commitment (UC) is conducted in a predict-then-optimize process: it starts by predicting the renewable energy source (RES) availability and system reserve requirements; given the predictions, the UC model is then optimized to determine the economic operation plans. In fact, predictions within the process are \textit{raw}. In other words, if the predictions are further tailored to assist UC in making the economic operation plans against realizations of the RES and reserve requirements, UC economics will benefit significantly. To this end, this paper presents a cost-oriented tailor of RES-and-reserve predictions for UC, deployed as an add-on to the predict-then-optimize process. The RES-and-reserve tailor is trained by solving a bi-level mixed-integer programming model: the upper level trains the tailor based on its induced operating cost; the lower level, given tailored predictions, mimics the system operation process and feeds the induced operating cost back to the upper level; finally, the upper level evaluates the training quality according to the fed-back cost. Through this training, the tailor learns to customize the raw predictions into cost-oriented predictions. Moreover, the tailor can be embedded into the existing predict-then-optimize process as an add-on, improving the UC economics. Lastly, the presented method is compared to traditional, binary-relaxing, neural network-based, stochastic, and robust methods.
\end{abstract}

\vspace{-3mm}
\section*{Nomenclature}
\addcontentsline{toc}{section}{Nomenclature}
\begin{spacing}{0.99}
\subsection*{Sets and Indexes}
\begin{IEEEdescription}[\IEEEusemathlabelsep \IEEEsetlabelwidth{$\mspace{45mu}$} \setlength{\IEEElabelindent}{0pt}]
\item[$\mathcal{B}/b$]
Set/index of branches.

\item[$\mathcal{E}/e$]
Set/index of iterations. The iteration limitation is $E$.

\item[$\mathcal{I}/i$]
Set/index of non-RES units, i.e., $\mathcal{I} = \mathcal{I}^{\text{ns}} \cup \mathcal{I}^{\text{qs}}$, where $\mathcal{I}^{\text{ns}}$ and $\mathcal{I}^{\text{qs}}$ are sets of non-quick-start and quick-start non-RESs.

\item[$\mathcal{J}/j$]
Set/index of RES units.

\item[$\mathcal{K}/k$]
Set/index of piecewise linear segments of generation curve.

\item[$\mathcal{S}/s$]
Set/index of training samples.

\item[$\mathcal{T}/t, t^{\prime}$]
Set/indexes of hours. The total number of hours is $T$.

\item[$\mathcal{T}^{\text{su/sd}}_{i}$]
Set defined as $\{{T}^{\text{su}}_{i},\cdots,T\}/\{{T}^{\text{sd}}_{i},\cdots,T\}$.

\item[$|\cdot|$]
The cardinality of a set.

\end{IEEEdescription}
\vspace{-4mm}
\subsection*{Decision Variables}
\begin{IEEEdescription}[\IEEEusemathlabelsep \IEEEsetlabelwidth{$\mspace{45mu}$} \setlength{\IEEElabelindent}{0pt}]
\item[$D_{it}$]
Shutdown status of unit $i$ at hour $t$.

\item[$I_{it}$]
On-off status of unit $i$ at hour $t$.

\item[$O_{it}$]
Non-spinning reserve (NR) commitment of unit $i$ at hour $t$.

\item[$P_{it}$]
Generation schedule of unit $i$ at hour $t$.

\item[$P_{itk}^{\text{sg}}$]
Generation schedule of unit $i$ in segment $k$ at hour $t$.

\item[$R_{it}^{\text{sr/nr}}$]
Spinning reserve (SR)/NR schedule of unit $i$ at hour $t$.

\item[$S$]
Slack variable.


\item[$U_{it}$]
Startup status of unit $i$ at hour $t$.

\item[$\mathcal{W}/\mathcal{R}$]
Tailor of RES power/reserve requirements.

\item[$W_{jt}$]
Generation schedule of RES $j$ at hour $t$.

\item[$\boldsymbol{x}, \boldsymbol{y}$]
Vectors of decision variables in the UC model, including binaries 
$\boldsymbol{x}=$
$\{\boldsymbol{U},$
$\boldsymbol{I},$
$\boldsymbol{O},$
$\boldsymbol{D}\}$ 
and continuous
$\boldsymbol{y} =
\{\boldsymbol{P},
\boldsymbol{P}^{\text{sg}},
\boldsymbol{R}^{\text{sr}},
\boldsymbol{R}^{\text{nr}},
\boldsymbol{W}\}$.

{\item[$\bar{\boldsymbol{x}}^{\text{ev}}$]
Enumerated solutions to binary variables.}

\item[$\boldsymbol{z}$]
Vector of decision variables in the re-dispatch (RD) model, including binaries $\{\boldsymbol{U}^{\text{RD}}, \boldsymbol{I}^{\text{RD}}, \boldsymbol{I}^{\text{RD,qs}}, \boldsymbol{D}^{\text{RD}}\}$ and continuous $\{ \boldsymbol{P}^{\text{RD}}, \boldsymbol{P}^{\text{RD,sg}}, \boldsymbol{W}^{\text{RD}}, \boldsymbol{S}\}$.

\item[${\cdot}^{\text{dv/gv}}$]
Duplicated/generated variables.

\item[${\cdot}^{\text{RD}}$]
Variables of RD problem.

\item[${\cdot}^{\star}$]
Indicating the optimal solution to a variable.

\end{IEEEdescription}

\vspace{-4mm}
\subsection*{Parameters}
\begin{IEEEdescription}[\IEEEusemathlabelsep \IEEEsetlabelwidth{$\mspace{45mu}$} \setlength{\IEEElabelindent}{0pt}]
\item[$B_b$]
Transmission capacity of branch $b$.

\item[$C_{ik}^{\text{sg}}$]
Generation price of unit $i$ in segment $k$.

\item[$C_{i}^{\text{su/nl}}$]
Startup/no-load price of unit $i$.

\item[$C^{\text{gs/ls/bs}}$]
Penalty price of slack variables.

\item[$\mathcal{P}^{\text{w/r}}$]
Feasible region of predictor $\mathcal{W}(\cdot)/\mathcal{R}(\cdot)$.

\item[$P^{\text{M/m}}_{i}$]
Maximum/minimum generation limit of unit $i$.

\item[$\bar{P}_{ik}^{\text{sg}}$]
Power limit of unit $i$ at segment $k$.

\item[$\bar{R}_{i}^{\text{sr/nr}}$]
SR/NR limit of unit $i$.

\item[$\hat{R}_{t}^{\text{sr/nr}}$]
Predicted SR/NR requirement at hour $t$, which together form the vector $\hat{\boldsymbol{r}}$.

\item[$R_{i}^{\text{su/sd}}$]
Startup/shutdown ramping capacity of unit $i$.

\item[$R_{i}^{\uparrow/\downarrow}$]
Upward/downward ramping capacity of unit $i$.

\item[$T_{i}^{\text{su/sd}}$]
Minimum on/off time requirement of unit $i$.

\item[$\hat{W}/\tilde{W}_{jt}$]
Predicted/actual available power of RES $j$ at hour $t$, forming the vector $\hat{\boldsymbol{w}}/\tilde{\boldsymbol{w}}$.

\item[$\hat{\boldsymbol{w}}/\hat{\boldsymbol{r}}$]
{\text{Raw prediction of RES/reserve requirements.}}

\item[$\hat{\boldsymbol{w}}^{\diamond}/\hat{\boldsymbol{r}}^{\diamond}$]
{\text{Tailored counterpart of raw prediction $\hat{\boldsymbol{w}}/\hat{\boldsymbol{r}}$.}}

\end{IEEEdescription}






















\end{spacing}

\section{Introduction}\label{sec01}
\subsection{Problem Statement}
\IEEEPARstart{T}{his} paper focuses on the unit commitment (UC) economics of power systems integrated with renewable energy sources (RES), where system operators are the top administrators \cite{Cost}. The daily UC is generally executed in a predict-then-optimize process, shown in Fig.~\ref{fig01}(a) and summarized as the following three steps:
\begin{figure}[tb]
	\centering
		\includegraphics[width=\columnwidth]{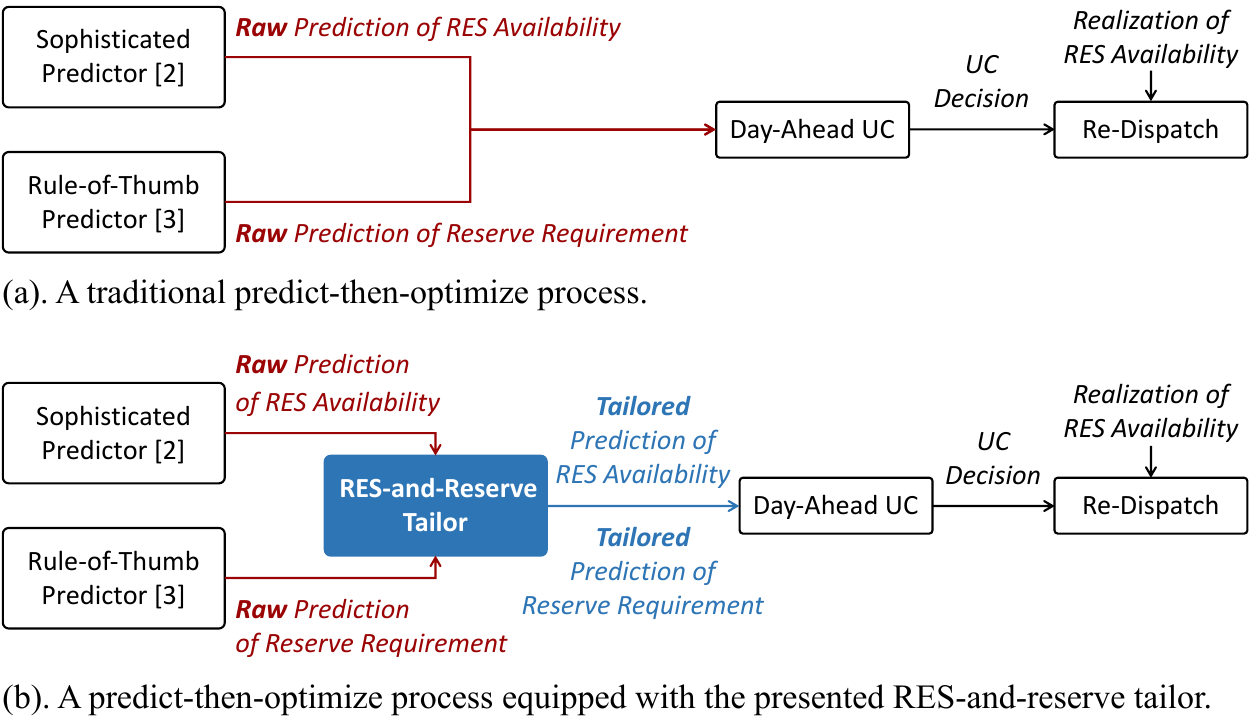}
		\vspace{-6mm}
	\caption{Predict-then-optimize processes with and without the add-on tailor.}\label{fig01}
	\vspace{-3mm}
\end{figure}

\begin{itemize}[noitemsep, topsep=0pt, parsep=0pt, partopsep=0pt]
\item
\textit{Step 1:} The RES availability \cite{Intro_RESPRE}, as well as the system spinning reserve (SR) and non-spinning reserve (NR) requirements \cite{ReservePre1}, are predicted in the day-ahead stage. These predictions are generated for serving various downstream applications, including but not limited to UC. Thus, they are regarded as \textit{raw predictions}\footnote{Raw predictions are typically provided by sophisticated \cite{Intro_RESPRE} or rule-of-thumb \cite{ReservePre1} predictors. Some systems make the raw predictions publicly accessible, enabling various entities to use these predictions for their specific applications.};
\item
\textit{Step 2:} Using the raw predictions as inputs, system operators solve a deterministic UC problem \cite{Intro_UC} to determine day-ahead operation plans involving startup/shutdown schedules and generation baselines;
\item
\textit{Step 3:} To evaluate the economics of a UC plan, a re-dispatch (RD) problem \cite{Intro_ED}, which is built upon the UC plan, is solved to calculate the eventual \textit{balancing cost} against the actual RES realizations as well as SR and NR deployments.
\end{itemize}

With respect to the day-ahead predictions and the afterward RES realizations, the startup and no-load costs from day-ahead UC and the balancing costs from RD constitute the \textit{actual operating cost}. This is referred to as \textit{UC economics} in this paper. Indeed, a UC plan with superior economics can enable a low actual operating cost.

Although the current predict-then-optimize process seems plausible, its UC economics actually suffers from the raw predictions for at least two reasons: \textit{i)} the RES predictors are generally trained by statistical accuracy criteria, such as mean absolute percentage error (MAPE). However, due to the inevitable prediction errors and the inherent nonlinearity of the UC-RD optimization process, \textit{a statistically more accurate prediction may NOT necessarily induce a lower actual operating cost} (as revealed in Section~\ref{asymmetric}); and \textit{ii)} the raw SR/NR requirements sized by traditional rule-of-thumb predictors could be overly redundant or insufficient, further worsening the UC economics.

Considering that the current raw predictors have been sophisticatedly developed and are utilized across various downstream applications, this paper targets to explore the following open question: \textit{How to leverage the predictors currently used by the existing predict-then-optimize practice and tailor their output, i.e., raw predictions, to improve the UC economics?}

{
\subsection{Literature Review}
To improve the UC economics against the inevitable inaccuracy of raw predictions, one approach is to augment the deterministic UC to its two-stage stochastic programming (TSP) \cite{Intro_SPUC} or two-stage robust optimization (TRO) \cite{RO1} counterpart. Differently, the \textit{closed-loop predict-and-optimize} idea, proposed by \cite{bengio} as an alternative approach, has gained growing attention. As illustrated in Fig.~\ref{fig02}, this idea emphasizes that \textit{predictors should support operators by prioritizing making good final decisions over pursuing immediate prediction accuracy, thus necessitating a closed-loop interaction between the prediction and optimization.} Indeed, this closed-loop idea has also been highlighted in the IEEE-CIS predict+optimize report \cite{challenge} and the INFORMS special series entitled ``Blending Predictive \& Prescriptive Methods'' \cite{informs}. Moreover, pioneering methodologies \cite{bertsimas, elmachtoub, ban} and software packages \cite{pkg1, pkg2, pkg3} have also been developed.
\begin{figure}[tb]
	\centering
		\includegraphics[width=\columnwidth]{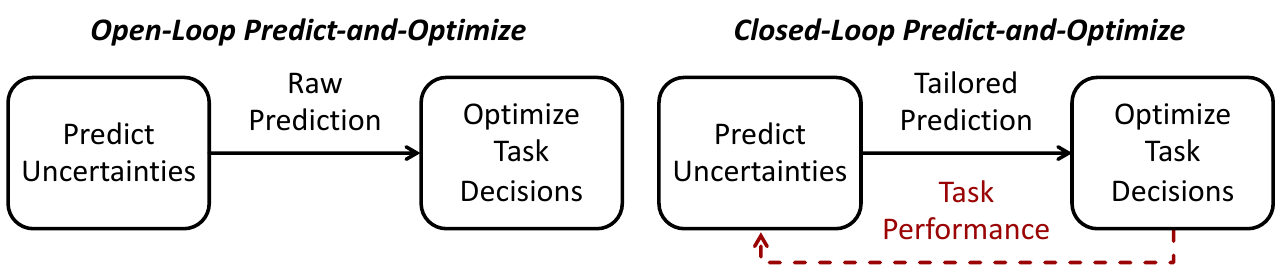}
		\vspace{-5mm}
	\caption{The ideas of open-loop and closed-loop predict-and-optimize.}\label{fig02}
	\vspace{-4mm}
\end{figure}

{
\begin{table*}[tb]
\vspace{-2mm}
	\caption{Review of Relevant Application References (Involving Closed-Loop Predict-and-Optimize Idea) in Power System Fields}\label{tab01}
	\centering
	\footnotesize
\begin{tabularx}{\textwidth}{c@{\extracolsep{\fill}}ccccccc}
\toprule
\multirow{2.5}{*}{Type}&\multirow{2.5}{*}{Ref.}&\multicolumn{2}{c}{Optimization Task}&\multirow{3}{*}{\shortstack[c]{Feedback Path between\\Prediction and Optimization}}&\multirow{2.5}{*}{Training Method} &\multirow{3}{*}{\shortstack[c]{Mathematical Model\\of ERM Problem}}\\
\cmidrule{3-4}
&                  &Task Content                 &Task Model       &                         &&\\
\midrule
\multirow{6.0}{*}{IPO}
&\cite{carriere}   &RES trading       &LP  &\multirow{6.5}{*}{\shortstack[c]{Merge the prediction step\\and optimization step as\\a single prescription step}}        &\multirow{2}{*}{Heuristic} &\multirow{2}{*}{-}\\
&\cite{stratigakos}&RES trading       &QP  & & &\\
\cmidrule{3-4} \cmidrule{6-7}
&\cite{sang1}      &Voltage regulation&SOCP& &Gradient-based updating &-\\
\cmidrule{3-4} \cmidrule{6-7}
&\cite{munoz1}     &RES trading       &LP  & &Optimally solve the ERM model      &Single-level LP\\
&\cite{changfei1}  &Reserve sizing    &LP  & &to get optimally-trained predictors    &Single-level MIP\\
\midrule

\multirow{21}{*}{CPO}

&\cite{yamin}  &UC&MIP &\multirow{3.5}{*}{\shortstack[c]{Use approximating error-and-cost\\function as loss function}}
                                                          &\multirow{2}{*}{Heuristic}&\multirow{2}{*}{-}\\
&\cite{guoli}  &RES trading         &LP&                          &                          &\\

\cmidrule{3-4} \cmidrule{6-7}

&\cite{yiwang1}&ED&LP&&Gradient-based updating&-\\

\cmidrule{2-7}
&\cite{mingyang1}&System planning &MIP&\multirow{3}{*}{\shortstack[c]{Identify scenarios that are\\helpful in reducing task cost}}                                   &\multirow{3}{*}{\shortstack[c]{Depend on the\\applied ML algorithm}}&\\
&\cite{zhangchao}&ED              &LP &      &                                                                    &-\\
&\cite{yurdakul} &UC              &MIP&      &                                                                    &\\

\cmidrule{2-7}

&\cite{donti,jiayu,yufan2}  &ED                      &LP &\multirow{6}{*}{\shortstack[c]{Solve the task optimization\\model with given predictions\\and then feed its cost back to\\the predictor training phase}}&\multirow{4}{*}{Gradient-based updating}&\multirow{4}{*}{-}\\
&\cite{wilder,mandi1,mandi2}&UC                      &MIP&&                                        &\\
&\cite{sang2}               &Storage arbitrage&MIP&&                                        &\\
&\cite{yiwang2}             &ED with storage  &MIP&&                                        &\\
\cmidrule{3-4} \cmidrule{6-7}

&\cite{yufan1}   &ED              &LP & &Contextual bandit-type &-\\

\cmidrule{2-7}
&\cite{mingyang2}   &Inertia management&NLP&\multirow{5.5}{*}{\shortstack[c]{Integrate\\predictor training and\\task optimization\\in the ERM model\\via explicit constraints}}&\multirow{6}{*}{\shortstack[c]{Optimally solve the\\ERM model to get\\optimally-trained predictors}}&Single-level NLP\\
&\cite{xianbang,wuhan}&UC              &MIP&                                     &                   &Single-level MIP\\
\cmidrule{3-4} \cmidrule{7-7}
&\cite{morales1}      &RES trading     &LP &                                     &                   &\multirow{2}{*}{\shortstack[c]{Bi-level LP}}\\
&\cite{munoz2, morales2, garcia, vladimir, viafora, changfei2}  
                      &ED              &LP &                                     &                   &\\

\cmidrule{3-4} \cmidrule{7-7}

&This paper          &UC \& RD         &MIP \& MIP&                              &                   &Bi-level MIP\\
\midrule
\multicolumn{7}{l}{QP: Quadratic programming; SOCP: Second-order cone programming; NLP: Non-linear programming.}\\
\bottomrule
\end{tabularx}
\vspace{-4mm}
\end{table*}

The closed-loop idea has been increasingly applied to power systems. Table~\ref{tab01} lists 30 power system-related references, compiled to the best of the authors' knowledge as of January 2024. These references are classified into \textit{``integrated predict-and-optimize (IPO)''} and \textit{``cost-oriented predict, then optimize (CPO).''} Both IPO \cite{carriere, stratigakos, sang1, munoz1, changfei1} and CPO \cite{yamin, guoli, yiwang1, mingyang1, zhangchao, yurdakul, donti,jiayu,yufan2,wilder,mandi1,mandi2, sang2, yiwang2, yufan1, mingyang2, xianbang,wuhan, morales1,munoz2, morales2, vladimir, garcia, viafora, changfei2} aim to improve the performance\textemdash generally in terms of the economics\textemdash of optimization tasks, while they adopt different philosophies:
\begin{itemize}[noitemsep, topsep=0pt, parsep=0pt, partopsep=0pt]
\item
IPO streamlines the predict-then-optimize process as a single prescription step, by using machine learning (ML) models such as neural networks (NN) to prescribe decisions from features directly. IPO-based frameworks are well-suited for tasks whose decision feasibility is easily reachable. However, they generally cannot assure decision feasibility when facing complex constraints;
\item
CPO follows the two-step predict-then-optimize practice but enhances the prediction step to induce superior decisions in the optimization step. CPO-based frameworks are distinguished by their custom-trained predictors, which can yield \textit{cost-oriented} predictions to enable superior decision-making in the optimization task. Although the tailored predictions may slightly compromise the statistical accuracy, they can boost the task performance greatly.
\end{itemize}

Early CPO frameworks \cite{yamin, guoli, yiwang1} approximate the relationship between prediction error and task cost via analytical functions, which are then used as training loss functions to replace traditional accuracy-oriented loss functions. Furthermore, references \cite{donti, jiayu, yufan2} train a novel cost-aware load/RES predictor for economic dispatch (ED). Their training takes the ED cost induced by predictions as the loss, and their predictor parameters are updated by the derivative of ED cost w.r.t. predictor parameters. Such a cost-aware predictor can generate predictions tailored to reduce the ED cost. However, references \cite{wilder, mandi1, mandi2} indicate that this approach is limited to linear programming (LP)-based operation tasks, primarily due to: \textit{i)} the need for gradient information in training; and \textit{ii)} the optimization task model is repeatedly solved per updating epoch, resulting in a time-consuming training. To adapt this approach for mixed-integer programming (MIP)-based tasks, great efforts have been made in \cite{wilder, mandi1, mandi2, sang2, yiwang2}. Nevertheless, the adaptability is achieved by either compromisingly relaxing/fixing integer variables \cite{wilder, mandi2, yiwang2} or ignoring the non-differentiable issue \cite{mandi1, sang2}.

In response to this limitation, references \cite{mingyang2, xianbang, wuhan} propose to use mathematical programming-based methods for training predictors. This training typically involves two steps: \textit{i)} constructing the empirical risk minimization (ERM) problem as an explicit mathematical model (e.g., MIP), which integrates the predictor training and task optimization in a closed-loop manner; and \textit{ii)} solving this ERM model to yield optimally trained predictors. As illustrated in \cite{mingyang2, xianbang, wuhan}, CPO frameworks based on this training inherit two key properties of mathematical programming-based methods: they can uncompromisingly handle the integer variables because the training does not depend on the gradient, and they exhibit consistent performance as the predictor is trained optimally.

Nevertheless, the ERM model in \cite{mingyang2, xianbang, wuhan} is being challenged recently. One major concern arises from their \textit{single-level} structure\textemdash a structure that cannot capture the chronological sequence between the prediction and optimization. This inability may incur an over-fitting predictor, especially when the training samples are limited, as emphasized in \cite{munoz2}. To overcome this issue, references \cite{morales1, munoz2, morales2, vladimir, garcia, viafora, changfei2} use \textit{bi-level} programming-based ERMs\textemdash training the predictor at the upper level while mimicking the optimization task at the lower level. Their numerical results demonstrate that CPO frameworks with the bi-level ERM model can improve task economics remarkably while mitigating the over-fitting issue.

Based on the above literature review and Table~\ref{tab01}, an important research gap can be revealed\textemdash the existing bi-level programming-based CPO frameworks \cite{munoz2, morales1, morales2, vladimir, garcia, viafora, changfei2} exclusively target LP tasks. This is largely because their training methods rely on the favorable properties of LP models. Typically, a lower-level LP subproblem can be equivalently replaced with its Karush–Kuhn–Tucker (KKT) conditions, thereby converting the original bi-level model into a tractable single-level form. However, it has been widely recognized that many practical operation tasks can only be reasonably formed as MIP problems, in which the KKT conditions cannot be applied directly. To this end, a broader adaptability of the bi-level programming-based CPO frameworks is necessary.

\subsection{Main Works and Contributions}
Motivated by the above-identified research gap, this work presents a bi-level MIP-based CPO framework, which trains a cost-oriented tailor of RES-and-reserve predictions for the MIP-based UC. The term \textit{``tailor''} refers to a custom-made predictor that can modify the raw predictions into cost-oriented predictions for improving the UC economics.

To train the tailor, the ERM problem is constructed in a bi-level MIP form. Its upper level trains the tailors for RES power as well as SR and NR requirements, in which the input/output is the raw/tailored prediction. The bi-level ERM employs \textit{two} MIP-based lower levels to mimic UC and RD operation processes. These lower levels utilize the upper-level predictions to calculate the actual operating cost, which is then fed back to the upper level for improving the tailoring quality.

Regarding the solving method for the ERM model, it begins by proving that optimal solutions of the original bi-level MIP (with two lower levels) can be obtained by solving a relatively tractable bi-level MIP (with only one lower level). The more tractable bi-level MIP is then solved by a cutting plane method, yielding the optimally trained tailors of RES, SR, and NR.

Finally, as shown in Fig.~\ref{fig01}(b), the tailors act as add-ons to the existing predict-then-optimize practice, positioned between the currently-used predictors and the UC model. The tailors can customize the raw RES-SR-NR predictions into their cost-oriented counterparts, thus improving the UC economics.

Our contributions are summarized as twofold:

\textit{1)} To improve the UC economics, a CPO framework including a bi-level MIP-based ERM model and its solution method is presented. Compared to the existing works, the proposed CPO framework is distinguished by the following aspects:
\begin{itemize}[noitemsep, topsep=0pt, parsep=0pt, partopsep=0pt]
{
\item
The presented bi-level MIP-based CPO framework can economically benefit UC and other MIP-based operation tasks in power systems. As compared, most gradient-based methods (e.g., \cite{donti}) and bi-level LP-based frameworks (e.g., \cite{morales1, munoz2, morales2, vladimir, garcia, viafora, changfei2}) are restricted to LP tasks or require compromises, e.g., relaxing binary variables;
}

\item
Recall that revealing the eventual UC economics necessitates solving the MIP-based RD model. To this end, the presented bi-level ERM model is designed to incorporate both UC and RD into the training as two lower-level subproblems. In comparison, the existing UC-related works (i.e., \cite{yamin, yurdakul, wilder, mandi1, mandi2, xianbang, wuhan}) neglect the RD task. Therefore, the presented framework is versatile for more general tasks;

\item
Theoretically, the presented mathematical programming-based training method can consistently approach a globally optimal predictor. Contrarily, the training in \cite{yamin, guoli, yiwang1, mingyang1, zhangchao, yurdakul, donti, jiayu, yufan2, wilder, mandi1, mandi2, sang2, yiwang2, yufan1} typically provides a locally-optimal predictor at random. This means that the presented framework can perform consistently, especially when the number of training samples is limited.
\end{itemize}

\textit{2)}
To illustrate the effectiveness of the presented CPO framework, comprehensive comparisons are conducted using real-world data. Specifically, the presented CPO framework is compared to the predict-then-optimize practice, bi-level relaxed MIP, NN-based CPO \cite{yiwang1}, TSP \cite{Intro_SPUC}, and TRO \cite{RO1}.}}

The rest of the paper is organized as follows: Section~\ref{sec02} details the necessary preliminaries; Section~\ref{sec03} expounds on the training of the tailor; Section~\ref{sec04} analyzes the numerical results; and Section~\ref{sec05} concludes this paper.

\section{Preliminaries}\label{sec02}
\subsection{Operation Models}
Taking the raw predictions $\hat{\boldsymbol{w}}$ and $\hat{\boldsymbol{r}}$ as inputs, the operators solve the MIP-based UC model to determine the day-ahead operation plans. The compact form of UC is shown as \eqref{UC}, and its detailed formulation is presented in Appendix~\ref{app_a}.
\begin{subequations}\label{UC}
\begin{align}
 \min_{\boldsymbol{x},\boldsymbol{y}}\,
          &\boldsymbol{b}^{\top}\boldsymbol{x} + \boldsymbol{c}^{\top}\boldsymbol{y}                  \label{UC:1}\\
 s.\,t.\, &\boldsymbol{x}, \boldsymbol{y} \in \mathcal{X}(\hat{\boldsymbol{w}}, \hat{\boldsymbol{r}});\label{UC:2}\\
          & \mathcal{X}(\hat{\boldsymbol{w}}, \hat{\boldsymbol{r}})=
                                         \left\{\boldsymbol{x}, \boldsymbol{y}\,\,
                                                \begin{array}{|l}
                                                      \boldsymbol{F}(\boldsymbol{x}, \boldsymbol{y}) = 0;  \\
                                                      \boldsymbol{G}(\boldsymbol{x}, \boldsymbol{y},
                                                                     \hat{\boldsymbol{w}}, \hat{\boldsymbol{r}}) \leq 0;
                                                \end{array}
                                         \right\};                                                       \label{UC:3}
\end{align}
\end{subequations}

Here, $\boldsymbol{F}(\cdot)$/$\boldsymbol{G}(\cdot)$ represents the set of UC equality/inequality constraints, forming the feasible region $\mathcal{X}$. With raw predictions $\hat{\boldsymbol{w}}$ and $\hat{\boldsymbol{r}}$ as inputs, UC is solved to deliver optimal solutions of startup, commitment, baseline generation, SR, and NR schedules \cite{shao1}.

Next, with respect to the operation plans $\{\boldsymbol{x}^{\star}, \boldsymbol{y}^{\star}\}$ from UC, the best balancing operations against actual RES realization $\tilde{\boldsymbol{w}}$ can be revealed by solving the MIP-based RD model. The compact form of the RD is shown as \eqref{ED}, and its detailed formulation is presented in Appendix~\ref{app_b}.
\begin{subequations}\label{ED}
\begin{align}
\min_{\boldsymbol{z}}\,
          & \boldsymbol{d}^{\top}\boldsymbol{z}                        \label{ED:1}\\
s.\,t.\,  & \boldsymbol{z}  \in \mathcal{Z}(\boldsymbol{x}^{\star},
                                            \boldsymbol{y}^{\star},
                                            \tilde{\boldsymbol{w}});   \label{ED:2}\\
          & \mathcal{Z}(\boldsymbol{x}^{\star},
                        \boldsymbol{y}^{\star},
                               \tilde{\boldsymbol{w}}) =
                           \left\{ \boldsymbol{z}\,\,
                                   \begin{array}{|l}
                                          \boldsymbol{M}(\boldsymbol{x}^{\star}, \boldsymbol{z}) = 0;\\
                                          \boldsymbol{N}(\boldsymbol{x}^{\star}, \boldsymbol{y}^{\star},
                                                        \boldsymbol{z},
                                                        \tilde{\boldsymbol{w}}) \leq 0
                                   \end{array}
                           \right\};                                  \label{ED:3}
\end{align}
\end{subequations}

Here, $\boldsymbol{M}(\cdot)$/$\boldsymbol{N}(\cdot)$ denotes the set of RD equality/inequality constraints, and $\mathcal{Z}$ is the feasible region. Solving RD provides optimal solutions to status switch of quick-start units, generation adjustments of all units, and slack penalties against the RES realization $\tilde{\boldsymbol{w}}$ \cite{zl2}.

\subsection{The Evaluation of UC Economics}\label{ucrdprocess}
The UC economics can be assessed in three steps:

\textit{i)} As the realization $\tilde{\boldsymbol{w}}$ remains unknown in the day-ahead stage, operators first conduct UC based on raw predictions $\hat{\boldsymbol{w}}$ and $\hat{\boldsymbol{r}}$. This step provides the \textit{actual UC cost} $\boldsymbol{b}^{\top}\boldsymbol{x}^{\star}$ consisting of startup and no-load costs;

\textit{ii)} After $\tilde{\boldsymbol{w}}$ is revealed, RD is executed based on the UC decisions $\boldsymbol{x}^{\star}$. This step yields the \textit{balancing cost} $\boldsymbol{d}^{\top}\boldsymbol{z}^{\star}$, including startup and no-load costs of quick-start units additionally committed in RD, generation costs of all units, and slack penalties;

\textit{iii)} Eventually, equation \eqref{systemcost} provides the \textit{best possible actual operating cost}, i.e., $c^{\text{act}}$, that is realistically achievable w.r.t. the UC solution and actual RES realization. This is referred to as the \textit{UC economics} in this paper.
\begin{align}\label{systemcost}
c^{\text{act}} = &\textstyle{\left.\sum\nolimits_{i \in \mathcal{I}}\sum\nolimits_{t \in \mathcal{T}}
     (  C^{\text{su}}_{i} U_{it}^{\star}
      + C^{\text{nl}}_{i} I_{it}^{\star}) \mspace{123mu} \right\}} \boldsymbol{b}^{\top}\boldsymbol{x}^{\star}  && \mspace{-75mu}   \notag \\
          &\left.\begin{array}{l}
\mspace{-45mu} + \sum\limits_{i \in \mathcal{I}}\sum\limits_{t\in\mathcal{T}}
             [(C_{i}^{\text{su}} U_{it}^{\text{RD}\star} + C_{i}^{\text{nl}} I_{it}^{\text{RD}\star})
             +\sum\limits_{ k\in\mathcal{K}}C_{ik}^{\text{sg}}P_{itk}^{\text{RD,sg}\star}] \\
\mspace{-45mu} + \sum\limits_{t \in \mathcal{T}}
             [  C^{\text{gs}}  S_{t}^{\text{gs}\star}
              + C^{\text{ls}}  S_{t}^{\text{ls}\star} + \sum\limits_{{b} \in \mathcal{B}}C^{\text{bs}}(S_{bt}^{+\star}+ S_{bt}^{-\star})]\\
                 \end{array} \mspace{-8mu} \right\}               \boldsymbol{d}^{\top}\boldsymbol{z}^{\star}  \notag \\
=         &\,\, \textstyle{  \boldsymbol{b}^{\top}\boldsymbol{x}^{\star}
                           + \boldsymbol{d}^{\top}\boldsymbol{z}^{\star}}
\end{align}


\subsection{The Asymmetric Relationship between UC Economics and RES Prediction Accuracy}\label{asymmetric}
Most RES predictors use symmetrical metrics to evaluate statistical prediction accuracy. However, the relationship between the prediction accuracy and UC economics is inherently asymmetrical. To show this point, multiple raw RES predictions $\hat{\boldsymbol{w}}$ are generated based on actual Belgian system data \cite{Intro_RESPRE} on 12/17/2020. These raw predictions are applied on an IEEE 14-bus system to simulate the existing predict-then-optimize process, as shown in Fig.~\ref{fig01}(a). Finally, the loss of UC economics is evaluated by \eqref{EconomicsLoss}, where $c^{\text{act,raw}}$ and $c^{\text{act,perf}}$ are actual operating costs induced by raw $\hat{\boldsymbol{w}}$ and perfect error-free prediction $\tilde{\boldsymbol{w}}$. The reserve requirements are set as raw.
\begin{equation}\label{EconomicsLoss}
\text{Loss of UC Economics} = \frac{{c^{\text{act,raw}} - c^{\text{act,perf}}}}{c^{\text{act,perf}}}\times 100\%
\end{equation}

\begin{figure}[t]
	\centering
			\includegraphics[width=\columnwidth]{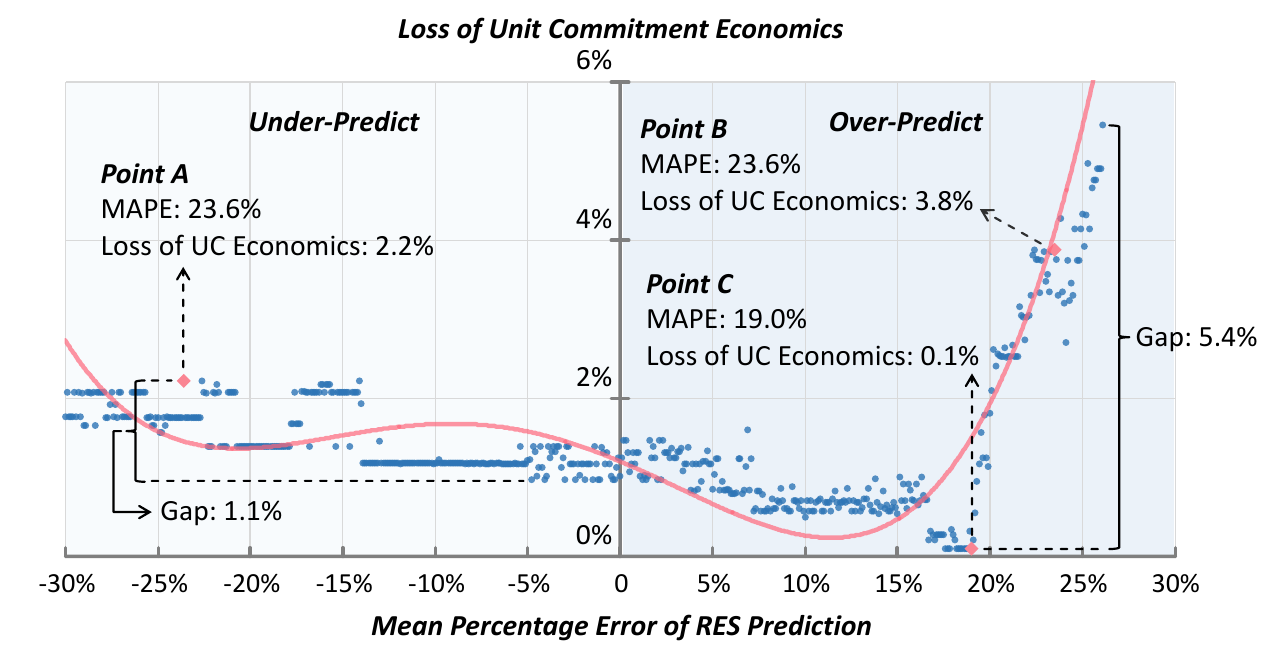}
		\vspace{-6mm}
	\caption{Impact of RES prediction error on UC economics.}\label{fig03}
\end{figure}

Fig.~\ref{fig03} shows the induced results of all raw prediction points, and the red line fits the relationship between MAPE and UC economics loss. The line clearly illustrates the asymmetry from a macro perspective: the loss caused by the over-predictions is fairly small in small-MAPE (5\% to 15\%) cases but becomes significant in large-MAPE (20\% to 30\%) cases, resulting in a 5.4\% width gap between the best and worst economic losses; in comparison, the loss induced by the under-predictions, the blue dots on the left half plane, is more stable with a narrower gap. Moreover, two cases are highlighted as follows:

\textit{Case 1: Raw predictions of the same MAPE may lead to different losses in UC economics.} Note that points \textit{A} and \textit{B} in Fig.~\ref{fig03} have the same MAPE (23.6\%) but different losses (2.2\% vs 3.8\%). This is because the under-prediction of point \textit{A} merely causes downward generation adjustment and RES curtailments in the RD; in comparison, the over-prediction of point \textit{B} necessitates deploying more expensive NR from quick-start units in the RD;

\textit{Case 2: Raw predictions of a larger MAPE may even lead to a smaller loss in UC economics.} Although the MAPE of point \textit{C} is noticeably high (19.0\%), it achieves the lowest loss (0.1\%) among all raw predictions. This is because each $\hat{\boldsymbol{w}}$ is associated with a just-enough reserve level $\hat{\boldsymbol{r}}^{\star}$ w.r.t. $\tilde{\boldsymbol{w}}$\textemdash a proper $\hat{\boldsymbol{r}}^{\star}$ depends on the difference between $\hat{\boldsymbol{w}}$ and $\tilde{\boldsymbol{w}}$ as well as other system specifications. Consequently, if $\{\hat{\boldsymbol{w}}, \hat{\boldsymbol{r}}^{\star}\}$ is taken as input to \eqref{UC}, the induced UC decisions can deploy $\hat{\boldsymbol{r}}^{\star}$ to enable the most effective utilization of available RES in RD operations, achieving the best UC economics.

Thus, predictors trained via symmetrical accuracy metrics, although suitable for general purposes, ignore the inherent asymmetry of the UC problem and the subsequent RD process. This fact inspires us to develop a RES-and-reserve tailor for the UC problem. The tailor can customize raw predictions into a cost-oriented counterpart by considering the subsequent RD result, thus boosting the UC economics.

\subsection{Major Flaws in Single-Level ERMs}\label{vsspo}
Model \eqref{SPO} is a typical single-level ERM, where $s$/$\mathcal{S}$ is the index/set of training samples; ${c}^{\text{perf}}_{s}$ is the objective value \eqref{UC:1} induced by the actual RES realization $\tilde{\boldsymbol{w}}_{s}$; $\mathcal{W}(\cdot)$ is the RES predictor to be trained, with $\mathcal{P}^{\text{w}}$ being its feasible region; $\boldsymbol{f}_{s}$ is the input feature. The ERM aims to train a predictor $\mathcal{W}(\cdot)$ that can induce the UC cost closer to ${c}^{\text{perf}}_{s}$ for individual samples.
\begin{subequations}\label{SPO}
\begin{flalign}
\min_{\mathcal{W}, \boldsymbol{x}_{s},\boldsymbol{y}_{s}}
          &\textstyle{\frac{1}{|\mathcal{S}|} \sum\limits_{s \in \mathcal{S}}}
                              ||  \boldsymbol{b}^{\top}\boldsymbol{x}_{s}
                                + \boldsymbol{c}^{\top}\boldsymbol{y}_{s}
                                - c^{\text{perf}}_{s}||_{1} \label{SPO:1}\\
s.\,t.\, & \mathcal{W}(\cdot) \in \mathcal{P}^{\text{w}};  \label{SPO:2}\\
         & \boldsymbol{x}_{s}, \boldsymbol{y}_{s} \in \mathcal{X}(\mathcal{W}(\boldsymbol{f}_{s}), \hat{\boldsymbol{r}}_{s}),\,\forall s \in \mathcal{S};                            \label{SPO:3}
\end{flalign}
\end{subequations}

The major flaw of \eqref{SPO} is that its in-sample UC solutions $\{\boldsymbol{x}_{s}^{\star}, \boldsymbol{y}_{s}^{\star}\}$ may use more expensive units, even if cheaper units are feasible, to pursue the minimum loss \eqref{SPO:1}. Such UC decisions violate the least-cost  principle\footnote{{{As explained in \cite{vladimir}, system operators are responsible for identifying the least-cost UC decision against the currently available predictions, known as the least-cost principle. Thus, operators generally use \eqref{UC:1} as the objective.}}} of UC models, thus being compromised.

An example with two samples, $s1$ and $s2$, can illustrate this flaw \cite{munoz2}. Both samples have the same features ($\boldsymbol{f}_{s1} = \boldsymbol{f}_{s2}$) but different RES realizations ($\tilde{\boldsymbol{w}}_{s1} \neq \tilde{\boldsymbol{w}}_{s2}$). With these two samples, the single-level model \eqref{SPO} \textit{cannot} ensure their in-sample UC decisions, $\{\boldsymbol{x}_{s1}^{\star}, \boldsymbol{y}_{s1}^{\star}\}$ and $\{\boldsymbol{x}_{s2}^{\star}, \boldsymbol{y}_{s2}^{\star}\}$, are the same. This fact can be explained by the following two steps:
\begin{itemize}[noitemsep, topsep=0pt, parsep=0pt, partopsep=0pt]
\item[\textit{i)}]
$\boldsymbol{f}_{s1} = \boldsymbol{f}_{s2}$ means $\mathcal{W}(\boldsymbol{f}_{s1}) = \mathcal{W}(\boldsymbol{f}_{s2})$. That is, their in-sample predictions are the same;

\item[\textit{ii)}]
$\tilde{\boldsymbol{w}}_{s1} \neq \tilde{\boldsymbol{w}}_{s2}$ will lead to $c^{\text{perf}}_{s1} \neq c^{\text{perf}}_{s2}$. Given this, \eqref{SPO} will derive different in-sample UC decisions $\{\boldsymbol{x}_{s1}^{\star}, \boldsymbol{y}_{s1}^{\star}\} \neq \{\boldsymbol{x}_{s2}^{\star}, \boldsymbol{y}_{s2}^{\star}\}$ to achieve the minimum training loss \eqref{SPO:1}. However, the UC model \eqref{UC}, which follows the least-cost principle, necessitates that the same predictions must lead to the same optimal UC decisions. Thus, it can be concluded that at least one of the two in-sample decisions is compromised by violating the least-cost principle.

\end{itemize}

From the perspective of mathematical modeling, this flaw arises because the predictor training and UC decision-making are packed into a single level without following the chronological order\textemdash they are mutually anticipatable\footnote{{Anticipativity and non-anticipativity are concepts in stochastic programming. They indicate that in a multi-stage model, the current stage can (anticipativity) or cannot (non-anticipativity) access and manipulate the information and decisions of future stages.}} in \eqref{SPO}. As a result, \eqref{SPO} will derive unequal UC solutions $\{\boldsymbol{x}_{s1}^{\star}, \boldsymbol{y}_{s1}^{\star}\} \neq \{\boldsymbol{x}_{s2}^{\star}, \boldsymbol{y}_{s2}^{\star}\}$ to chase a lower training loss \eqref{SPO:1}. In other words, in-sample UC decisions of \eqref{SPO} are made to adapt the in-sample cost $c^{\text{perf}}$ instead of minimizing operating cost, thus failing to reflect the least-cost principle-based decision-making of operators. As a result, predictors trained via the single-level ERM \eqref{SPO}, although presenting good in-sample losses, could fail to generalize well (i.e., over-fitting) because the in-sample UC decisions are compromised.

In comparison, constructing ERM as a bi-level MIP can layer the prediction and optimization problems in a proper chronological order, thus satisfying the non-anticipativity naturally. In this way, the in-sample UC decisions always minimize the operating cost \eqref{UC:1}, thus strictly obeying the least-cost principle and effectively mitigating the over-fitting issue.

\section{RES-and-Reserve Tailors based on\\ Bi-Level MIP}\label{sec03}
\subsection{Constructing the ERM Problem}
Using a pre-specified dataset $\mathcal{S}$ (Section~\ref{sampleselect} will detail the dataset construction) as well as raw predictions $\hat{\boldsymbol{w}}_{s}$ and $\hat{\boldsymbol{r}}_{s}$, the bi-level ERM problem is constructed as in \eqref{ERM}. The upper level \eqref{ERM:1}-\eqref{ERM:3} trains tailors $\mathcal{W}$ and $\mathcal{R}$ for RES, SR, and NR, taking actual operating costs \eqref{systemcost} as the loss function \eqref{ERM:1}. The loss function maximizes the weighted sum of UC economics for all samples in $\mathcal{S}$. Constraint \eqref{ERM:2} limits the tailor parameters within their feasible regions $\mathcal{P}^{\text{w}}$ and $\mathcal{P}^{\text{r}}$; \eqref{ERM:3} defines $\mathcal{W}$ and $\mathcal{R}$, which respectively take raw predictions $\hat{\boldsymbol{w}}_{s}$ and $\hat{\boldsymbol{r}}_{s}$ as inputs. The two lower levels, \eqref{ERM:4} and \eqref{ERM:5}, model UC \eqref{UC} and RD \eqref{ED}, respectively. Solving \eqref{ERM} can provide cost-oriented tailors, $\mathcal{W}^{\star}$ and $\mathcal{R}^{\star}$, that are trained optimally.
\begin{subequations}\label{ERM}
\begin{flalign}
&\textit{Upper Level (Training Tailors of Raw Predictions):}
                                                         \mspace{-100mu}&                                     \notag\\
&\min_{\mathcal{W}, \mathcal{R}} \textstyle{\frac{1}{|\mathcal{S}|} \sum\limits_{s \in \mathcal{S}}} {\boldsymbol{b}^{\top} \boldsymbol{x}_{s}} + {\boldsymbol{d}^{\top} \boldsymbol{z}_{s}}
                                                        \mspace{-100mu}&                               \label{ERM:1}\\
&s.\,t.\,\mathcal{W}(\cdot) \in \mathcal{P}^{\text{w}},\, \mathcal{R}(\cdot) \in \mathcal{P}^{\text{r}};
                                                         \mspace{-100mu}&                                \label{ERM:2}\\
&\mspace{30mu}\hat{\boldsymbol{w}}^{\diamond}_{s}=\mathcal{W}(\hat{\boldsymbol{w}}_{s}),\, \hat{\boldsymbol{r}}^{\diamond}_{s}=\mathcal{R}(\hat{\boldsymbol{r}}_{s}),                 & \forall s \in \mathcal{S};               \label{ERM:3}\\
&\mspace{30mu}\textit{Lower Level (Day-Ahead Unit Commitment):}
                                                         \mspace{-100mu}&                                      \notag\\
&\mspace{30mu}\boldsymbol{x}_{s}, \boldsymbol{y}_{s} \in \underset{\boldsymbol{x}_{s}, \boldsymbol{y}_{s} \in \mathcal{X}(\hat{\boldsymbol{w}}^{\diamond}_{s}, \hat{\boldsymbol{r}}^{\diamond}_{s})}{\arg \min}\boldsymbol{b}^{\top}\boldsymbol{x}_{s}+\boldsymbol{c}^{\top} \boldsymbol{y}_{s},          \mspace{-100mu}&  \forall s \in \mathcal{S};    \label{ERM:4}\\
&\mspace{30mu}\textit{Lower Level (Re-Dispatch):}        \mspace{-100mu}&                                   \notag\\
&\mspace{30mu}\boldsymbol{z}_{s} \in \underset{\boldsymbol{z}_{s} \in \mathcal{Z}(\boldsymbol{x}_{s}, \boldsymbol{y}_{s}, \tilde{\boldsymbol{w}}_{s})}{\arg\min} \boldsymbol{d}^{\top} \boldsymbol{z}_{s},
                                                         \mspace{-100mu}&     \forall s \in \mathcal{S}; \label{ERM:5}
\end{flalign}
\end{subequations}

The bi-level ERM \eqref{ERM} features that: \textit{i)} its bi-level structure and the $\arg \min$ operators of \eqref{ERM:4} and \eqref{ERM:5} enable the in-sample decisions to obey the least-cost principle; and \textit{ii)} it can incorporate two MIP-based tasks into the training.

\begin{figure}[tb]
	\centering
		\includegraphics[width=\columnwidth]{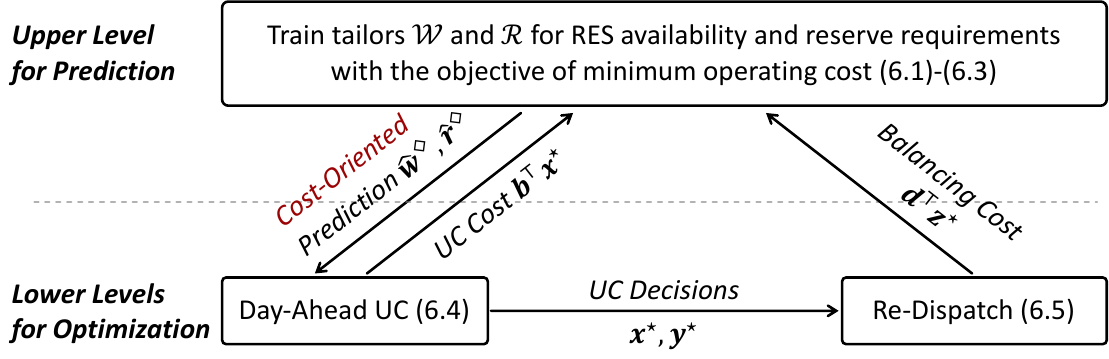}
		\vspace{-4mm}
	\caption{The closed-loop predict-and-optimize mechanism in the ERM \eqref{ERM}.}\label{fig04}
\end{figure}

Fig.~\ref{fig04} expounds how the bi-level ERM \eqref{ERM} achieves the closed-loop predict-and-optimize idea. The upper level \eqref{ERM:1} generates the tailored predictions $\hat{\boldsymbol{w}}^{\diamond}$ and $\hat{\boldsymbol{r}}^{\diamond}$, but cannot directly assess their induced UC economics. Thus, $\hat{\boldsymbol{w}}^{\diamond}$ and $\hat{\boldsymbol{r}}^{\diamond}$ are first passed to the lower-level UC \eqref{ERM:4}, forming the UC feasible region $\mathcal{X}$. Solving \eqref{ERM:4} provides the actual UC cost $\boldsymbol{b}^{\top} \boldsymbol{x}^{\star}$ to the upper level as well as the UC decisions $\{\boldsymbol{x}^{\star}, \boldsymbol{y}^{\star}\}$ to the lower-level RD \eqref{ERM:5}. Taking $\{\boldsymbol{x}^{\star}, \boldsymbol{y}^{\star}\}$ as inputs, \eqref{ERM:5} is solved to further provide the balancing cost $\boldsymbol{d}^{\top} \boldsymbol{z}^{\star}$ to the upper level. Finally, the upper level can assess the prediction quality via the actual operating cost $\boldsymbol{b}^{\top} \boldsymbol{x}^{\star} + \boldsymbol{d}^{\top} \boldsymbol{z}^{\star}$, and tune the tailor parameters accordingly. This process repeats until the tailors $\mathcal{W}$ and $\mathcal{R}$ are trained optimally. Since tailors $\mathcal{W}^{\star}$ and $\mathcal{R}^{\star}$ are trained in favor of the UC economics, they are cost-oriented customization for the UC problem.

\subsection{On Selecting Specific Models for Tailor}
The tailors act as add-ons in the existing predict-then-optimize process, in which they output prediction vectors $\hat{\boldsymbol{w}}^{\diamond}_{s}$ and $\hat{\boldsymbol{r}}^{\diamond}_{s}$ for the day-ahead UC. Given that these tailors are tasked with predicting continuous values, their training could be better treated as a \textit{regression problem}.

In theory, the ERM model \eqref{ERM} can accommodate various \textit{regression forms} by formulating them as explicit constraints \eqref{ERM:3} to be handled by the solvers. That is, \eqref{ERM:3} can take various forms, such as linear, polynomial, logarithmic, and parametric ReLU. Among these forms, this paper selects the affine linear form \eqref{form}.\begin{subequations}\label{form}\begin{flalign}
&\mathcal{W}(\hat{\boldsymbol{w}}_s) = \boldsymbol{m} \circ \hat{\boldsymbol{w}}_s,\,\boldsymbol{m}\in \mathbb{R}_{+}^{(|\mathcal{T}|  \times |\mathcal{J}|) \times 1};\label{form:1}\\
&\mathcal{R}(\hat{\boldsymbol{r}}_s) = \boldsymbol{n} \circ \hat{\boldsymbol{r}}_s,\,\boldsymbol{n}\in \mathbb{R}_{+}^{2|\mathcal{T}| \times 1};\label{form:2}
\end{flalign}
\end{subequations}

In \eqref{form}, $\circ$ is the element-wise product, and the tailor parameters are constrained to nonnegative real numbers. With this, the essence of the tailoring is to linearly scale raw predictions in favor of UC economics. The training is basically to solve the ERM \eqref{ERM} for determining the optimal vectors $\boldsymbol{m}$ and $\boldsymbol{n}$.

Notably, the topic of selecting a proper regression form has attracted increasing attention within the IPO and CPO fields. This interest mainly stems from the fact that IPO and CPO methods are designed for a specific task (e.g., UC), necessitating their forms to align with the task properties.

Specifically, the prescription model of IPO methods directly maps features to the ultimate operation decisions. Consequently, complex forms, such as NN \cite{sang1}, are typically used to capture the highly complicated and non-linear relationship. On the other hand, the prediction model of CPO methods first transforms the input feature into cost-oriented prediction, which is subsequently used as input to operation models for optimizing ultimate decisions. In this case, the form of the prediction model mainly depends on the relationship between the input feature and the cost-oriented prediction.

In the presented CPO method, the high-quality raw predictions sourced from sophisticated predictors, rather than basic meteorological data like temperature, are used as features. For such a case, \cite{munoz1, xianbang}, and \cite{munoz2, morales2, garcia} have demonstrated that the raw prediction generally requires only minor adjustments to render its cost-oriented counterpart. This is why the paper uses ``tailor'' instead of ``predictor'' to describe $\mathcal{W}$ and $\mathcal{R}$. Thus, the straightforward linear form \eqref{form} would be suitable for the presented method due to its balanced effectiveness, efficiency, interpretability, and numerical stability. 

\subsection{Solving the ERM Problem}\label{Solving}
The primary computational challenge in solving the ERM problem \eqref{ERM} arises from the binary variables in the two lower-level problems, i.e., \eqref{ERM:4} and \eqref{ERM:5}. If no binary variables exist, \eqref{ERM:4} and \eqref{ERM:5} can be equivalently converted to their KKT conditions, enabling \eqref{ERM} to be solved as an MIP problem.

The presented method takes two steps to handle the above challenge: \textit{i)} proving that the optimal solutions of \eqref{ERM} can be obtained by solving a more tractable bi-level MIP with only one lower-level subproblem; and \textit{ii)} solving this more tractable bi-level MIP with only one lower-level subproblem through a cutting plane-based method.

\subsubsection{Converting (6) with Two Lower-Level MIP Subproblems to (8) with One Lower-Level MIP Subproblem} The bi-level ERM \eqref{ERM} includes two MIP subproblems in the lower level, introducing non-negligible computational complexity. 

To solve \eqref{ERM}, Proposition~\ref{P1} is proposed to convert \eqref{ERM} into a more tractable bi-level MIP \eqref{ERMMM} with one MIP in the lower level. Solving \eqref{ERMMM} provides the optimal solution of the original ERM \eqref{ERM}, as proven in Proposition~\ref{P1}.
\begin{proposition}\label{P1}
The optimal solution of the bi-level MIP \eqref{ERMMM} is also optimal for the original bi-level MIP \eqref{ERM}.
\begin{subequations}\label{ERMMM}
\begin{flalign}
{\min_{\mathcal{W},
        \mathcal{R},
        \boldsymbol{z}_{s}}}
        &\textstyle{\frac{1}{|\mathcal{S}|} \sum\limits_{s \in \mathcal{S}}}
         \boldsymbol{b}^{\top} \boldsymbol{x}_{s}
                  + \boldsymbol{d}^{\top} \boldsymbol{z}_{s} \mspace{-45mu}& \label{ERMMM:1}\\
s.\,t.\,&\mathcal{W}(\cdot) \in \mathcal{P}^{\text{w}},\,
         \mathcal{R}(\cdot) \in \mathcal{P}^{\text{r}};    \mspace{-45mu}& \label{ERMMM:2}\\
        &\hat{\boldsymbol{w}}^{\diamond}_{s} = \mathcal{W}(\hat{\boldsymbol{w}}_{s}),\,
            \hat{\boldsymbol{r}}^{\diamond}_{s} = \mathcal{R}(\hat{\boldsymbol{r}}_{s}),
                                                             \mspace{-45mu}&\forall s \in \mathcal{S}; \label{ERMMM:3}\\
        &\boldsymbol{z}_{s} \in \mathcal{Z}(\boldsymbol{x}_{s}, \boldsymbol{y}_{s},
                          \tilde{\boldsymbol{w}}_{s}),       \mspace{-45mu}&\forall s \in \mathcal{S}; \label{ERMMM:4}\\
        &\boldsymbol{x}_{s}, \boldsymbol{y}_{s} \in
                \underset{\boldsymbol{x}_{s}, \boldsymbol{y}_{s} \in \mathcal{X}(\hat{\boldsymbol{w}}^{\diamond}_{s},
                           \hat{\boldsymbol{r}}^{\diamond}_{s})}{\arg \min} \boldsymbol{b}^{\top} \boldsymbol{x}_{s}
                 + \boldsymbol{c}^{\top} \boldsymbol{y}_{s}, \mspace{-45mu}& \forall s \in \mathcal{S}; \label{ERMMM:5}
\end{flalign}
\end{subequations}
\end{proposition}

\begin{IEEEproof}
To begin with, model \eqref{ERMM} is introduced to bridge \eqref{ERM} and \eqref{ERMMM}, which is constructed by adding \eqref{ERMM:4} to \eqref{ERM}, or equivalently, by adding \eqref{ERMM:6} to \eqref{ERMMM}.
\begin{subequations}\label{ERMM}
\begin{flalign}
 \min_{\mathcal{W},\mathcal{R},\boldsymbol{z}_{s}}
         &\textstyle{\frac{1}{|\mathcal{S}|} \sum\limits_{s \in \mathcal{S}}}
                         \boldsymbol{b}^{\top} \boldsymbol{x}_{s} + \boldsymbol{d}^{\top} \boldsymbol{z}_{s}
                         \mspace{-45mu}&                             \label{ERMM:1}\\
 s.\,t.\,&\mathcal{W}(\cdot) \in \mathcal{P}^{\text{w}},\, \mathcal{R}(\cdot) \in \mathcal{P}^{\text{r}};
                         \mspace{-45mu}&                             \label{ERMM:2}\\
         &\hat{\boldsymbol{w}}^{\diamond}_{s}=\mathcal{W}(\hat{\boldsymbol{w}}_{s}),\,
          \hat{\boldsymbol{r}}^{\diamond}_{s}=\mathcal{R}(\hat{\boldsymbol{r}}_{s}),
                         \mspace{-45mu}& \forall s \in \mathcal{S};  \label{ERMM:3}\\
         &\boldsymbol{z}_{s} \in \mathcal{Z}(\boldsymbol{x}_{s},\boldsymbol{y}_{s},\tilde{\boldsymbol{w}}_{s}),
                         \mspace{-45mu}& \forall s \in \mathcal{S};  \label{ERMM:4}\\
         &\boldsymbol{x}_{s},\boldsymbol{y}_{s} \in \underset{\boldsymbol{x}_{s},\boldsymbol{y}_{s} \in
                                   \mathcal{X}(\hat{\boldsymbol{w}}^{\diamond}_{s},
                                               \hat{\boldsymbol{r}}^{\diamond}_{s})}{\arg\min}
                                \boldsymbol{b}^{\top} \boldsymbol{x}_{s}
                              + \boldsymbol{c}^{\top} \boldsymbol{y}_{s},
                         \mspace{-45mu}& \forall s \in \mathcal{S};  \label{ERMM:5}\\
         &\boldsymbol{z}_{s} \in \underset{\boldsymbol{z}_{s}\in \mathcal{Z}(\boldsymbol{x}_{s},\boldsymbol{y}_{s},
                                 \tilde{\boldsymbol{w}}_{s})}{\arg\min} \boldsymbol{d}^{\top} \boldsymbol{z}_{s},
                         \mspace{-45mu}& \forall s \in \mathcal{S};  \label{ERMM:6}
\end{flalign}
\end{subequations}

First, we claim that \eqref{ERM} and \eqref{ERMM} have the same optimal solutions. This can be proved by showing that they have identical objective functions as well as feasible regions. It can be directly observed that their objective functions, \eqref{ERM:1} and \eqref{ERMM:1}, are the same. To show that \eqref{ERM} and \eqref{ERMM} have the same feasible region, we point out the following two facts: \textit{i)} a feasible solution of \eqref{ERM} satisfies all constraints of \eqref{ERMM} and thus is feasible to \eqref{ERMM}; and \textit{ii)} a feasible solution of \eqref{ERMM} also satisfies all constraints of \eqref{ERM} and thus is feasible to \eqref{ERM}. Thus, it is proven that \eqref{ERM} and \eqref{ERMM} have the same optimal solutions.

Second, we claim that the optimal solution of \eqref{ERMMM}, denoted as $\Omega^{\star}=\{ \boldsymbol{x}_{s}^{\star}, \boldsymbol{y}_{s}^{\star}, \boldsymbol{z}_{s}^{\star}\}$, is also optimal for \eqref{ERMM}. Note that since \eqref{ERMMM} is equivalent to \eqref{ERMM:1}-\eqref{ERMM:5}, $\Omega^{\star}$ is optimal for \eqref{ERMM:1}-\eqref{ERMM:5}. Thus, we only need to show that $\Omega^{\star}$ also satisfies \eqref{ERMM:6}, or equivalently, that $\boldsymbol{z}^{\star}_{s}$ satisfies \eqref{ERMM:6}. This can be achieved by three steps: \textit{i)} construct model \eqref{UCED} by removing \eqref{ERMM:2}-\eqref{ERMM:5} from \eqref{ERMM} and fixing the UC decisions at their optima $\boldsymbol{x}_{s}^{\star}$ and $\boldsymbol{y}_{s}^{\star}$. Note that \eqref{UCED} is equivalent to \eqref{ERMM:6}; \textit{ii)} recall that $\Omega^{\star}$ is an optimal solution of \eqref{ERMMM}, thus, its $\boldsymbol{z}^{\star}_{s}$ must be optimal for \eqref{UCED}; \textit{iii)} by combining steps \textit{i)} and \textit{ii)}, it is direct to conclude that $\boldsymbol{z}^{\star}_{s}$ satisfies \eqref{ERMM:6}.

As we have proved that $\Omega^{\star}$ is optimal for \eqref{ERMM:1}-\eqref{ERMM:5} and its $\boldsymbol{z}^{\star}_{s}$ also satisfies \eqref{ERMM:6}, it can be claimed that the optimal solution $\Omega^{\star}$ of \eqref{ERMMM} must be optimal for \eqref{ERMM}.
\begin{subequations}\label{UCED}
\begin{align}
 {\min_{\boldsymbol{z}_{s}}} &\textstyle{\frac{1}{|\mathcal{S}|}  \sum\limits_{s \in \mathcal{S}}}
                            \boldsymbol{b}^{\top} \boldsymbol{x}_{s}^{\star}
                           + \boldsymbol{d}^{\top} \boldsymbol{z}_{s}\\
s.\,t.\,                   &\textstyle{\boldsymbol{z}_{s} \in 
                                      \underset{\boldsymbol{z}_{s} \in
                                       \mathcal{Z}(\boldsymbol{x}^{\star}_{s},
                                                   \boldsymbol{y}^{\star}_{s},
                                                   \tilde{\boldsymbol{w}}_{s})}{\arg \min}
                             \boldsymbol{d}^{\top} \boldsymbol{z}_{s}}, \, \forall s \in \mathcal{S};
\end{align}
\end{subequations}

Finally, by combining the first claim (i.e., \eqref{ERM} and \eqref{ERMM} have the same optimal solutions) and the second claim (i.e., the optimal solution of \eqref{ERMMM} is also optimal for \eqref{ERMM}), it is direct to conclude that the optimal solution of \eqref{ERMMM} is optimal for \eqref{ERM}.
\end{IEEEproof}


\subsubsection{Solving (8) by C\&CG-based Algorithm}
\begin{algorithm}[t]
	\SetAlgoLined     
	\KwIn{Desired optimality gap $\epsilon^{\star}$,
	      initial tailors $\mathcal{W}_{e=0}$ and $\mathcal{R}_{e=0}$,
	      initial lower/upper bound $LB$/$UB$.}
	\KwOut{Optimally trained tailors $\mathcal{W}^{\star}$ and $\mathcal{R}^{\star}$.}
	\While{$e \leq E$}{
        Solve \textbf{SP1} and \textbf{SP2} for each sample $s \in \mathcal{S}$\;
		$UB \gets \min \{UB, \sum\limits_{s \in \mathcal{S}} \frac{c_{s}^{\text{act}\star}}{|\mathcal{S}|} + \lambda^{\text{w}} \lVert \mathcal{W}\rVert_{1} - \lambda^{\text{r}} \lVert\mathcal{R}\rVert_{1} \}$\;
		$\epsilon \gets 100\% \times \frac{UB - LB}{UB}$\;
		\eIf{$\epsilon \leq \epsilon^{\star}$}{
			$\mathcal{W}^{\star} \gets \mathcal{W}_{e}$, $\mathcal{R}^{\star} \gets \mathcal{R}_{e}$\;
			\textbf{break}\;
			}{
			${\bar{\boldsymbol{x}}_{s, e}^{\text{ev}}}$ $\gets$ optimal solution of decision $\boldsymbol{x}_{s}$ in \textbf{SP2}\;
			Add new $\{\boldsymbol{x}^{\text{gv}}_{s, e}, \boldsymbol{y}^{\text{gv}}_{s, e}\}$ and cuts \eqref{MP:6}-\eqref{MP:10} to \textbf{MP}\;
			Solve \textbf{MP}\;
			$e \gets e + 1, LB \gets$ optimal objective value of \textbf{MP}\;
			$\mathcal{W}_{e}$ and $\mathcal{R}_{e}$     $\gets$ optimal $\mathcal{W}$ and $\mathcal{R}$ from \textbf{MP}\;
		}}
	\caption{A C\&CG algorithm for solving \eqref{ERMMM}}\label{alg1}
\end{algorithm}

The bi-level MIP \eqref{ERMMM} can be solved by either Benders decomposition \cite{bard} or column-and-constraint generation (C\&CG). Since Benders decomposition may converge slowly, a C\&CG-based algorithm is selected, which iterates between the master problem \textbf{MP} and two subproblems \textbf{SP1} and \textbf{SP2}, as summarized in Algorithm~\ref{alg1}.

The first subproblem \textbf{SP1} is presented as in \eqref{SP1}. It is a duplication of the UC problem \eqref{ERMMM:5}, which takes the incumbent tailors $\mathcal{W}_{e}$ and $\mathcal{R}_{e}$ as inputs. Solving \textbf{SP1} can reveal the least-cost UC decision induced by the incumbent tailors, as well as its associated cost $c^{\text{SP1}\star}_{s}$.
\begin{subequations}\label{SP1}
\begin{flalign}
\textbf{SP1:}\, c^{\text{SP1}\star}_{s}=
        &\min_{\boldsymbol{x}_{s}, \boldsymbol{y}_{s}}
               \boldsymbol{b}^{\top}\boldsymbol{x}_{s}
               + \boldsymbol{c}^{\top}\boldsymbol{y}_{s}                                 \label{SP1:1}\\
s.\,t.\,&\hat{\boldsymbol{w}}^{\diamond}_{s}=\mathcal{W}_{e}(\hat{\boldsymbol{w}}_{s}),\,
         \hat{\boldsymbol{r}}^{\diamond}_{s}=\mathcal{R}_{e}(\hat{\boldsymbol{r}}_{s});  \label{SP1:2}\\
        &\boldsymbol{F}(\boldsymbol{x}_{s}, \boldsymbol{y}_{s}) = 0,\,\boldsymbol{G}(\boldsymbol{x}_{s},
                        \boldsymbol{y}_{s},
                        \hat{\boldsymbol{w}}^{\diamond}_{s},
                        \hat{\boldsymbol{r}}^{\diamond}_{s}) \leq 0;                     \label{SP1:3}
\end{flalign}
\end{subequations}

Note that multiple optimal binaries $\boldsymbol{x}_{s}^{\star}$ for \textbf{SP1} may exist, and they could lead to different objective values of \eqref{ERMMM:1}. Thus, the second subproblem \textbf{SP2}, formed as in \eqref{SP2} with $c^{\text{SP1}\star}_{s}$ as input, is further solved to select a solution among multiple $\boldsymbol{x}_{s}^{\star}$ of \textbf{SP1} in favor of the objective \eqref{ERMMM:1}. This is referred to as \textit{optimistic bi-level programming} \cite{ref25}.
\begin{subequations}\label{SP2}
\begin{flalign}
\textbf{SP2:}\,
c_{s}^{\text{act}\star} =
        &\min_{\boldsymbol{x}_{s},
               \boldsymbol{y}_{s}, \boldsymbol{z}_{s}}
         \boldsymbol{b}^{\top} \boldsymbol{x}_{s}
          + \boldsymbol{d}^{\top} \boldsymbol{z}_{s}                                      \label{SP2:1}\\
s.\,t.\,&\hat{\boldsymbol{w}}^{\diamond}_{s} = \mathcal{W}_{e}(\hat{\boldsymbol{w}}_{s}),\,
         \hat{\boldsymbol{r}}^{\diamond}_{s} = \mathcal{R}_{e}(\hat{\boldsymbol{r}}_{s}); \label{SP2:2}\\
        &\boldsymbol{b}^{\top} \boldsymbol{x}_{s}
         + \boldsymbol{c}^{\top} \boldsymbol{y}_{s} \leq c^{\text{SP1}\star}_{s};       \label{SP2:3}\\
        &\boldsymbol{F}(\boldsymbol{x}_{s},\boldsymbol{y}_{s}) = 0,\,
         \boldsymbol{G}(\boldsymbol{x}_{s},\boldsymbol{y}_{s},
         \hat{\boldsymbol{w}}^{\diamond}_{s},
         \hat{\boldsymbol{r}}^{\diamond}_{s}) \leq 0 ;                                    \label{SP2:4}\\
        &\boldsymbol{M}(\boldsymbol{x}_{s}, \boldsymbol{z}_{s}) = 0,\,
         \boldsymbol{N}(\boldsymbol{x}_{s}, \boldsymbol{y}_{s}, \boldsymbol{z}_{s},
                        \tilde{\boldsymbol{w}}_{s}) \leq 0;                               \label{SP2:5}
\end{flalign}
\end{subequations}

Given the binary decisions $\boldsymbol{x}_{s}^{\star}$ identified by \textbf{SP2}, {the master problem \textbf{MP} is formulated as \eqref{MP}, which is composed of the \textit{high-point problem} \cite{ref25} \eqref{MP:1}-\eqref{MP:5} of \eqref{ERMMM} and the cutting planes \eqref{MP:6}-\eqref{MP:10}.} Three variants of variables $\{\boldsymbol{x}_{s}, \boldsymbol{y}_{s}\}$ are used in \textbf{MP}, including duplicated variables $\{\boldsymbol{x}^{\text{dv}}_{s},\boldsymbol{y}^{\text{dv}}_{s}\}$, enumerated binary variables {$\bar{\boldsymbol{x}}^{\text{ev}}_{s,e}$}, and generated continuous variables $\boldsymbol{y}^{\text{gv}}_{s,e}$. The duplicated variables serve as the proxy of the original variables, but with a larger feasible region; the enumerated and generated variables and their associated cutting planes \eqref{MP:6}-\eqref{MP:10} serve to cut the feasible region of duplicated variables gradually. \textbf{MP} is built in three parts:
\begin{itemize}[noitemsep, topsep=0pt, parsep=0pt, partopsep=0pt]
\item
The first part \eqref{MP:1}-\eqref{MP:4} is based on the duplicated variables. {The objective function \eqref{MP:1} includes two regularization terms, $\lambda^{\text{w}} \lVert \mathcal{W} \rVert_{1}$ and $\lambda^{\text{r}} \lVert \mathcal{R} \rVert_{1}$, to mitigate over-predictions and insufficient reserves via hyper-parameters $\lambda^{\text{w}}$ and $\lambda^{\text{r}}$.} The regularization terms could be removed when the number of training samples is sufficient. Constraints \eqref{MP:2}-\eqref{MP:4} are the original upper-level constraints \eqref{ERMMM:2}-\eqref{ERMMM:4};

\item
The second part \eqref{MP:5} is a duplication of UC constraints to ensure the feasibility of the duplicated variables;

\item
The third part \eqref{MP:6}-\eqref{MP:10} acts as the \textit{optimality cuts}. {This part is built by three steps: \textit{i)} duplicating the MIP-based UC subproblem \eqref{ERMMM:5} for each sample $s$ in each iteration $e$; \textit{ii)} for each duplicated UC subproblem, fixing its binary variables as $\bar{\boldsymbol{x}}^{\text{ev}}_{s,e}$ that are given by \textbf{SP2}, which makes these MIPs degenerate to LPs; and \textit{iii)} substituting these LPs by their KKT conditions \eqref{MP:6}-\eqref{MP:9}.} Here, $\mathcal{L}(\cdot)$ denotes the Lagrangian function; $\boldsymbol{\mu}$ and $\boldsymbol{\nu}$ are dual variables of equality and inequality constraints. The KKT conditions \eqref{MP:6}-\eqref{MP:9} ensure the optimality of the generated variables w.r.t. {$\bar{\boldsymbol{x}}^{\text{ev}}_{s,e}$}. Finally, \eqref{MP:10} links the duplicated variables with their optimality cuts.
\end{itemize}
\begin{subequations}\label{MP}
\begin{flalign}
& \textbf{MP:}\, \min_{\Xi} \textstyle{\frac{1}{|\mathcal{S}|} \sum\nolimits_{s \in \mathcal{S}}}
  (  \boldsymbol{b}^{\top}\boldsymbol{x}^{\text{dv}}_{s}
   + \boldsymbol{d}^{\top}\boldsymbol{z}_{s})
   + \lambda^{\text{w}} \lVert \mathcal{W}\rVert_{1} - \lambda^{\text{r}} \lVert\mathcal{R}\rVert_{1}
                                                \mspace{-245mu}&                                                \notag\\
& \text{where $\Xi = \{ \mathcal{W}(\cdot), \mathcal{R}(\cdot),
\hat{\boldsymbol{w}}^{\diamond}_{s},
\hat{\boldsymbol{r}}^{\diamond}_{s},
\boldsymbol{z}_{s}, 
\boldsymbol{x}^{\text{dv}}_{s},
\boldsymbol{y}^{\text{dv}}_{s},
\boldsymbol{y}^{\text{gv}}_{s} \}$}             \mspace{-245mu}&                                          \label{MP:1}\\
&\textit{Subject to:}                           \mspace{-245mu}&                                                \notag\\
&\mspace{10mu}\textit{Original Constraints:}    \mspace{-245mu}&                                                \notag\\
&\mspace{10mu}\mathcal{W}(\cdot) \in \mathcal{P}^{\text{w}},\,\mathcal{R} (\cdot) \in \mathcal{P}^{\text{r}};
                                                \mspace{-245mu}&                                          \label{MP:2}\\
&\mspace{10mu}\hat{\boldsymbol{w}}^{\diamond}_{s} = \mathcal{W}(\hat{\boldsymbol{w}}_{s}),\,\hat{\boldsymbol{r}}^{\diamond}_{s} = \mathcal{R}(\hat{\boldsymbol{r}}_{s}), \mspace{-245mu}&               \forall s \in \mathcal{S}; \label{MP:3}\\
&\mspace{10mu}\boldsymbol{M}(\boldsymbol{x}^{\text{dv}}_{s}, \boldsymbol{z}_{s}) = 0,\,
\boldsymbol{N}(\boldsymbol{x}^{\text{dv}}_{s},\boldsymbol{y}^{\text{dv}}_{s}, \boldsymbol{z}_{s}, \tilde{\boldsymbol{w}}_{s}) \leq 0,                                    \mspace{-245mu}&                \forall s \in \mathcal{S};\label{MP:4}\\
&\mspace{10mu}\textit{Duplication of UC Constraints:}
                                                 \mspace{-245mu}&                                           \notag\\
&\mspace{10mu}\boldsymbol{F}(\boldsymbol{x}_{s}^{\text{dv}}, \boldsymbol{y}_{s}^{\text{dv}}) = 0,\, \boldsymbol{G}(\boldsymbol{x}^{\text{dv}}_{s}, \boldsymbol{y}_{s}^{\text{dv}}, \hat{\boldsymbol{w}}_{s}^{\diamond}, \hat{\boldsymbol{r}}_{s}^{\diamond}) \leq 0,                              \mspace{-245mu}&                \forall s \in \mathcal{S};\label{MP:5}\\
&\mspace{10mu}\textit{Stationarity of KKT (Optimality Cuts):}
                                                \mspace{-245mu}&                                                \notag\\
&\mspace{10mu}\nabla\mathcal{L}(\boldsymbol{y}_{s,e}^{\text{gv}}, \boldsymbol{\mu}_{s,e}, \boldsymbol{\nu}_{s,e}) = 0,
                                                \mspace{-245mu}&                \forall s \in \mathcal{S};\label{MP:6}\\
&\mspace{10mu}\textit{Primal and Dual Feasibilities of KKT
         (Optimality Cuts):}                    \mspace{-245mu}&                                                \notag\\
&\mspace{10mu}\boldsymbol{F}({\bar{\boldsymbol{x}}_{s,e}^{\text{ev}}}, \boldsymbol{y}_{s,e}^{\text{gv}}) = 0,
                                                \mspace{-245mu}&  \forall s\in\mathcal{S},e\in\mathcal{E};\label{MP:7}\\
&\mspace{10mu}\boldsymbol{G}({\bar{\boldsymbol{x}}_{s,e}^{\text{ev}}}, \boldsymbol{y}_{s,e}^{\text{gv}},
      \hat{\boldsymbol{w}}_{s}^{\diamond},\hat{\boldsymbol{r}}_{s}^{\diamond}) \leq 0, \, \boldsymbol{\nu}_{s,e} \geq 0,
                                                \mspace{-245mu}& \forall s\in\mathcal{S},e\in\mathcal{E};\label{MP:8}\\
&\mspace{10mu}\textit{Complementary Slackness of KKT (Optimality Cuts):}
                                                \mspace{-245mu}&                                               \notag \\
&\mspace{10mu}\boldsymbol{\nu}_{s,e}\bot \boldsymbol{G}({\bar{\boldsymbol{x}}_{s,e}^{\text{ev}}}, \boldsymbol{y}_{s,e}^{\text{gv}}, \hat{\boldsymbol{w}}_{s}^{\diamond}, \hat{\boldsymbol{r}}_{s}^{\diamond}),
                                                \mspace{-245mu}&  \forall s\in\mathcal{S},e\in\mathcal{E};\label{MP:9}\\
&\mspace{10mu}\textit{Objective Cuts (Optimality Cuts):}
                                                \mspace{-245mu}&                                               \notag\\
&\mspace{10mu}\boldsymbol{b}^{\top} \boldsymbol{x}_{s}^{\text{dv}} + \boldsymbol{c}^{\top} \boldsymbol{y}_{s}^{\text{dv}} \leq \boldsymbol{b}^{\top} {\bar{\boldsymbol{x}}_{s,e}^{\text{ev}}} + \boldsymbol{c}^{\top} \boldsymbol{y}_{s,e}^{\text{gv}},
                                                \mspace{-245mu}&\forall s\in \mathcal{S},e\in\mathcal{E};\label{MP:10}
\end{flalign}
\end{subequations}

By linearizing the complementary slackness \eqref{MP:9} via the Big-M method, \textbf{MP} can be converted into an equivalent MIP problem and directly solved by commercial MIP solvers.


\subsection{Embedding RES-and-Reserve Tailor into UC as Add-On}
Recall that the trained cost-oriented tailors, $\mathcal{W}^{\star}$ and $\mathcal{R}^{\star}$, are linear functions \eqref{form}. Hence, they can be embedded into the deterministic UC \eqref{UC} as linear constraints \eqref{PUC:C}, resulting in a prescriptive UC \eqref{PUC}.
\begin{subequations}\label{PUC}
\begin{align}
\min_{\boldsymbol{x},\boldsymbol{y}, \hat{\boldsymbol{w}}^{\diamond}, \hat{\boldsymbol{r}}^{\diamond}}\,
        & \boldsymbol{b}^{\top}\boldsymbol{x}
                              + \boldsymbol{c}^{\top}\boldsymbol{y}                        \label{PUC:A}\\
s.\,t.\,& \boldsymbol{x},\boldsymbol{y} \in \mathcal{X}(\hat{\boldsymbol{w}}^{\diamond},
                                                        \hat{\boldsymbol{r}}^{\diamond});  \label{PUC:B}\\
        & \hat{\boldsymbol{w}}^{\diamond} = \mathcal{W}^{\star}(\hat{\boldsymbol{w}}),\,
          \hat{\boldsymbol{r}}^{\diamond} = \mathcal{R}^{\star}(\hat{\boldsymbol{r}});     \label{PUC:C}
\end{align}
\end{subequations}

The prescriptive UC features that: \textit{i)} without introducing extra computation burden, the new linear constraints \eqref{PUC:C} can tailor the raw predictions into their cost-oriented counterpart in favor of the UC economics; and \textit{ii)} the objective \eqref{PUC:A} strictly obeys the least-cost principle, thus \eqref{PUC} is compatible with the predict-then-optimize practice, as illustrated in Fig.~\ref{fig01}.

Finally, system operators can solve the prescriptive UC model \eqref{PUC}, instead of the deterministic UC model \eqref{UC}, to boost the UC economics. This is equivalent to using $\mathcal{W}^{\star}$ and $\mathcal{R}^{\star}$ to tailor the raw predictions into cost-oriented predictions and solve the corresponding deterministic UC \eqref{UC}.

\section{Case Studies}\label{sec04}
\subsection{Experimental Setting}

\subsubsection{A Weekly-Rolling Training for Tailors}\label{sampleselect}
To leverage compatibility with the operators' practice and ensure solution quality of \eqref{PUC}, the tailors are re-trained on a weekly basis to update \eqref{PUC} based on the past $\textit{NT}$ days. Specifically, to study the dispatch week of days $(D+\text{1})$ to $(D+\text{7})$, the past $\textit{NT}$ days (i.e., days $(D-NT)$ to $(D-\text{1}))$ are selected to build the training dataset $\mathcal{S}$. Then, the ERM \eqref{ERM} is solved to provide the cost-oriented tailors, $\mathcal{W}^{\star}$ and $\mathcal{R}^{\star}$, for the targeted week. At the end of this period, the above process is repeated to update the tailors for the next week, i.e., days $(D+\text{8})$ to $(D+\text{14})$.
\begin{table}[H]
\vspace{-2mm}
	\caption{Selective Weeks for Out-of-Sample Testing}\label{tab02}
	\centering
	\footnotesize
\begin{tabularx}{\columnwidth}{@{\extracolsep{\fill}}cccc}
\toprule
    Quarter 1, 2020     & Quarter 2, 2020    & Quarter 3, 2020    & Quarter 4, 2020 \\
\midrule
   Feb 04 - Feb 10    & Jun 28 - Jul 04  & Aug 20 - Aug 26  & Dec 12 - Dec 18\\
\bottomrule
\end{tabularx}
\end{table}

\subsubsection{Data Setup}
Data, including raw predictions and actual realizations, are collected from the data platform of a Belgian system \cite{Intro_RESPRE}. Based on the data, four weeks in 2020, as listed in Table~\ref{tab02}, are selected, representing peak RES availabilities in the four quarters. The data are scaled to build the RES and load data for the testing systems. The slack penalty prices are all set as \$2,000/MWh. The rule-of-thumb reserve requirement is set according to the CAISO practice \cite{ReservePre1}: the system reserve requirement, represented as a proportion $\alpha$ of the total system load, is met by non-RES units, of which at least 50\% is SR.

\subsubsection{Methods to be Compared}
\begin{figure*}[tb]
	\centering
	\includegraphics[width=\textwidth]{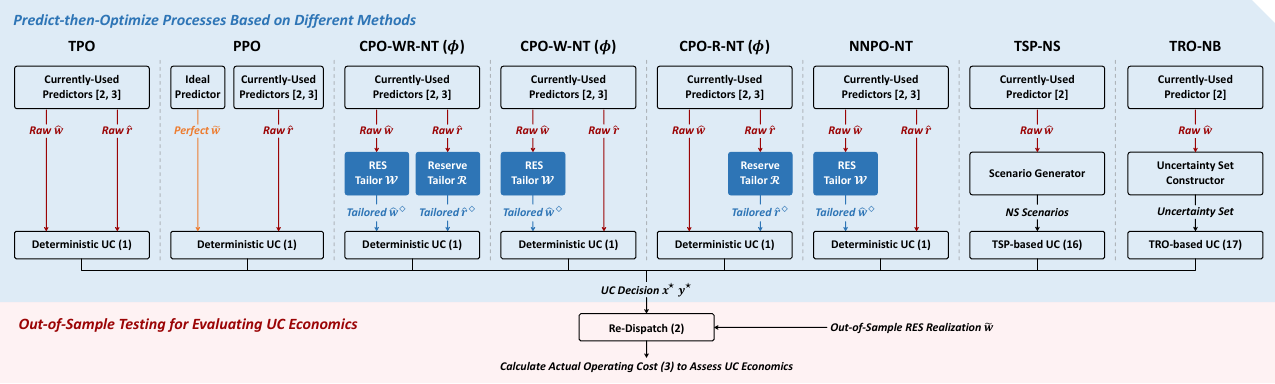}
	\vspace{-5mm}
	\caption{The methods to be analyzed and the comparison scheme.}\label{fig05}
\end{figure*}

Various methods are compared and analyzed to discuss the advantages and disadvantages of the presented framework. Fig.~\ref{fig05} shows the overall comparison strategy, which follows the evaluation process introduced in Section~\ref{ucrdprocess}. Note that \textit{NT} denotes the number of training samples. The methods are introduced as follows:
\begin{itemize}[noitemsep, topsep=0pt, parsep=0pt, partopsep=0pt]
\item
\textit{TPO:} Traditional predict-then-optimize. It uses raw predictions as inputs to \eqref{UC}, describing the current practice of system operators;

\item
\textit{PPO:} Perfect TPO. It uses error-free RES predictions and raw reserve requirements, representing the best UC economics an operator can achieve;

\item
\textit{CPO-WR-NT:} TPO enhanced with the cost-oriented RES and reserve tailors. This is the presented method;

\item
\textit{CPO-W-NT:} TPO enhanced with the cost-oriented RES tailor only\textemdash a variant of CPO-WR-NT. It tailors RES predictions but keeps the raw reserve requirements;

\item
\textit{CPO-R-NT:} TPO enhanced with the cost-oriented reserve tailor only\textemdash a variant of CPO-WR-NT. It keeps raw RES predictions but tailors the raw reserve requirements;

\item
\textit{CPO-NT-$\phi$:}
The same approach as CPO-NT, except that the binary requirement of lower-level operation models, \eqref{ERM:4} and \eqref{ERM:5}, are relaxed as in \eqref{LPUC} and \eqref{LPRD}. Parameter $\phi$ is pre-determined and falls within [0, 1], representing the ``relaxing degree'': $\phi$ = 1 means a 100\% relaxing degree\textemdash the lower-level models become LP; $\phi$ = 0 means no (0\%) relaxation. A larger $\phi$ indicates a higher relaxing degree of the binary requirement, i.e., the binary variables are allowed to take more fractional values between 0 and 1. Appendix~\ref{app_c} details this method;

\item
\textit{NNPO-NT \cite{yiwang1}:} TPO enhanced with an NN-based cost-oriented RES tailor. NNPO-NT is trained using the same samples as those of CPO-W-NT. Appendix~\ref{app_d} details this method;

\item
\textit{TSP-NS \cite{Intro_SPUC}:} Within the TPO framework, the deterministic UC model \eqref{UC} is augmented to its TSP counterpart. The first stage models UC, while the second stage simulates RD (w.r.t. the first-stage UC variables) in \textit{NS} net-load scenarios. Thus, the induced UC decisions can adapt to net-load realizations in the \textit{NS} scenarios;

\item
\textit{TRO-NB \cite{RO1}:} Within the TPO framework, the deterministic UC model \eqref{UC} is augmented to its TRO counterpart. The first stage models UC, while the second stage simulates RD (w.r.t. the first-stage UC variables) in the worst-case net-load scenario. \textit{NB} denotes the budget to handle the worst case.

\end{itemize}

The threshold on the optimality gap of all methods is set as 1\%. For CPO-WR-NT, the hyper-parameters, $\lambda^{\text{w}}$ and $\lambda^{\text{r}}$, are set as $\text{10}^{\text{6}}$. This is because the magnitude of the ERM objective \eqref{ERMMM:1} is between  $\text{10}^{\text{5}}$ and  $\text{10}^{\text{6}}$, so a hyper-parameter of $\text{10}^{\text{5}}$ would pose a sufficient, but not overly dominating, effect on the ERM. All cases are solved by Gurobi 9.1, called by YALMIP on MATLAB, on a 3.5 GHz PC.

\subsection{CPO vs TPO}
Based on an IEEE 14-bus system (with $\alpha$ = 50\%, 1,700MW non-RES capacity, 5 RES units, and 400MW RES capacity), this subsection compares CPO and TPO via economics improvement (EI) \eqref{EI} and value of tailoring (VoT) \eqref{VoT}.
\begin{align}
&\text{EI}  = \frac{c^{\text{act,TPO}}-c^{\text{act,CPO}}}{c^{\text{act,TPO}}}  \times 100\% \label{EI} \\
&\text{VoT} = \frac{c^{\text{act,TPO}}-c^{\text{act,CPO}}}{c^{\text{act,TPO}} - c^{\text{act,PPO}}}\label{VoT}
\end{align}

\subsubsection{Results}
\begin{figure}[H]
	\centering
		\includegraphics[width=\columnwidth]{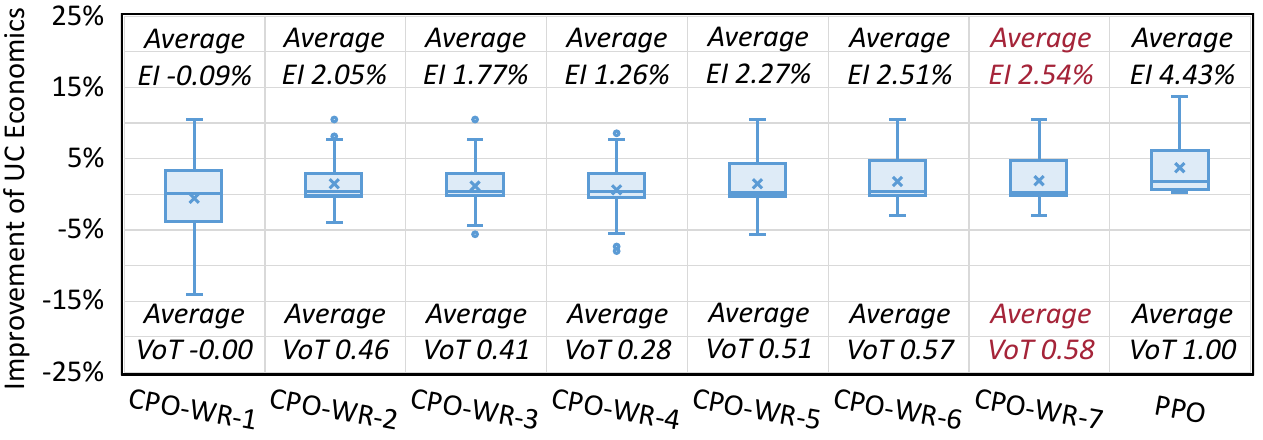}
		\vspace{-5mm}
	\caption{Out-of-sample UC economics of seven CPO-WR-NTs.}\label{fig06}
	\vspace{-2mm}
\end{figure}

Fig.~\ref{fig06} shows that six out of the seven CPO-WR-NTs, with \textit{NT} ranging from 2 to 7, outperform TPO economically. Notably, CPO-WR-7 renders the best performance\textemdash an average EI of 2.54\%. As the ideal PPO with the 1.00 VoT only achieves a 4.43\% EI, the 2.54\% EI of CPO-WR-7 can be regarded as non-insignificant. In fact, the 2.54\% EI corresponds to a considerable monthly cost reduction of \$171,360 compared to TPO.

\begin{figure}[H]
	\centering
		\includegraphics[width=\columnwidth]{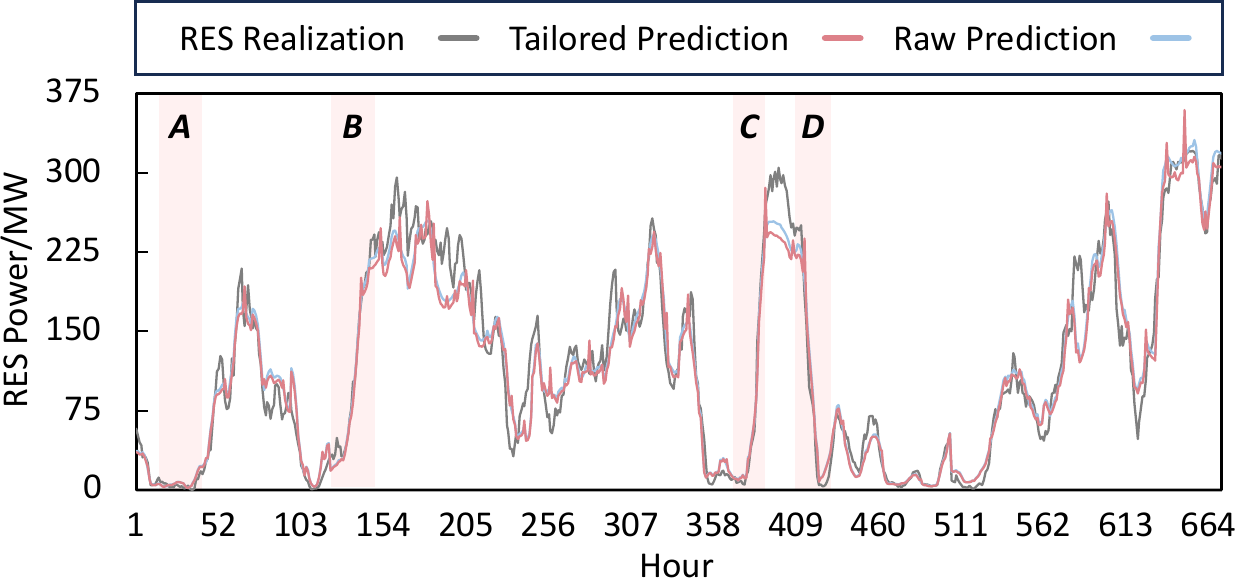}
		\vspace{-5mm}
	\caption{Raw (from TPO) and tailored (from CPO-WR-7) RES predictions.}\label{fig07}
\end{figure}

\begin{table}[H]
\vspace{-2mm}
	\caption{Breakdown of Out-of-Sample Costs of TPO and CPO-WR-7}\label{tab03}
	\centering
	\footnotesize
\begin{tabularx}{\columnwidth}{@{\extracolsep{\fill}} lcccccc}
\toprule
        &\multicolumn{2}{c}{Actual UC Cost/10$^{\text{3}}$\$}
        &$\mspace{-22mu}$&\multicolumn{2}{c}{Balancing Cost/10$^{\text{3}}$\$}
        &\multirow{4}{*}{$c^{\text{act}}$/10$^{\text{3}}$\$}\\
\cmidrule{2-3} \cmidrule{5-6}
        &\multirow{2}{*}{Startup}&\multirow{2}{*}{No-Load}&$\mspace{-22mu}$&Startup     &\multirow{2}{*}{Gen.}&\\
        &                        &                        &$\mspace{-22mu}$&plus No-Load&                     &\\

\midrule
TPO     &1.1                     &2.1                     &$\mspace{-22mu}$&0.2         & 237.9               &241.3\\
CPO-WR-7&1.0                     &2.2                     &$\mspace{-22mu}$&0.2         & 231.9               &235.3\\
\bottomrule
\end{tabularx}
\end{table}

Table~\ref{tab03} breaks down the actual operating costs of TPO and CPO-WR-7. The main difference between TPO and CPO-WR-7 stems from the generation cost (the $\text{5}^{\text{th}}$ column) in RD. It is also noteworthy that CPO-WR-7 has a lower startup cost and a higher no-load cost in UC.

Table~\ref{tab04} compares six CPO-NTs. Note that both CPO-WR-NT and CPO-R-NT render lower SR and NR than the raw schedules of TPO, among which SR is always higher than NR. Regarding the training, both the computational time and number of iterations for solving the ERM problem increase with \textit{NT}. Moreover, it is noteworthy that four CPO-NTs cause negative EI.

Table~\ref{tab05} evaluates RES prediction accuracy via five metrics, including mean absolute error (MAE), root mean square error (RMSE), MAPE, mean over-prediction percentage error (MOPE), and mean under-prediction percentage error (MUPE). Compared to TPO, CPO-WR-7 has a better MAPE but worse MAE and RMSE. Moreover, CPO-WR-7 has a better MOPE but a worse MUPE, indicating that CPO-WR-7 tends to conservatively predict RES power (i.e., vector $\hat{\boldsymbol{w}}^{\diamond}$ is statistically smaller than vector $\hat{\boldsymbol{w}}$).

Finally, Fig.~\ref{fig07} compares raw and tailored RES predictions (including 672 hours) of the four weeks, where two observations are noteworthy. First, the predictions after tailoring are not significantly different from raw predictions. This is because raw predictions are already of high accuracy, especially when RES is scarce (e.g., area \textit{A}) or the trend is stable (e.g., areas \textit{B}, \textit{C}, and \textit{D}), so there is no need to tailor them dramatically. Indeed, most elements of the vector $\boldsymbol{m}$ in \eqref{form} are between 0.9 and 1.1. Second, 86.9\% of tailored predictions are lower than their raw predictions, and the tailoring causes a larger deviation from the realizations.

\begin{table}[t]
	\caption{Comparison of Different CPO-NTs}\label{tab04}
	\centering
	\footnotesize
\begin{tabularx}{\columnwidth}{@{\extracolsep{\fill}}lccccr}
\toprule
         & \multicolumn{2}{c}{Scheduled Reserve}  &\multicolumn{2}{c}{Training}  &\multirow{2}{*}{EI/\%}         \\
\cmidrule{2-3} \cmidrule{4-5}
         &SR/MW       & NR/MW                   &Time/s & Iterations        &         \\
\midrule
TPO     &113.2       &113.2                    &-       &-                    &0.00\\
\midrule
CPO-WR-1&110.8       &75.6                     &19      &2                    &-0.09\\
CPO-WR-7&92.9        &65.4                     &786     &6                    &2.54 \\
\midrule
CPO-W-1 &113.2       &113.2                    &21      &2                    &-0.22\\
CPO-W-7 &113.2       &113.2                    &750     &3                    &-0.02 \\
\midrule
CPO-R-1 &98.8        &85.8                     &13      &2                    &-0.48\\
CPO-R-7 &88.8        &63.3                     &1,857    &4                    &1.77\\
\bottomrule
\end{tabularx}
\vspace{-1mm}
\end{table}

\begin{table}[t]
	\caption{Out-of-Sample Prediction Accuracy of TPO and CPO-WR-7}\label{tab05}
	\centering
	\footnotesize
\begin{tabularx}{\columnwidth}{@{\extracolsep{\fill}}lccccc}
\toprule
                                         &MAE            &RMSE           &MAPE           &MOPE        &MUPE\\
\midrule
Raw      $\hat{\boldsymbol{w}}$           &\textbf{15.0}MW&\textbf{21.0}MW&36.4\%         &29.6\%      &\textbf{6.8}\%\\
Tailored $\hat{\boldsymbol{w}}^{\diamond}$&16.5MW         &23.2MW         &\textbf{35.9}\%&\textbf{27.7}\%&8.2\%\\
\bottomrule
\end{tabularx}
\vspace{-1mm}
\end{table}

\subsubsection{Discussions} The above results have indicated that CPO-WR-NT can outperform TPO on the UC economics. Still, the reasons behind this conclusion deserve further discussion.

Table~\ref{tab03} shows that the cost difference mainly stems from the generation. This is because CPO-WR-7 tends to use units with cheaper generation costs (e.g., large-capacity units). Although these units have higher no-load costs, CPO-WR-7 keeps them online long enough by setting predictions $\hat{\boldsymbol{w}}^{\diamond}$ and $\hat{\boldsymbol{r}}^{\diamond}$ properly. As a result, CPO-WR-7 uses fewer units in RD and operates them most efficiently.

According to Table~\ref{tab04}, the results of CPO-WR-NT and CPO-R-NT reveal that a proper reserve deployment should have a higher SR and a lower NR, instead of the rule-of-thumb 50\%-50\% as in TPO. This is because the bi-level ERM \eqref{ERM} enables CPO-NT to learn that triggering quick-start units in RD is expensive. Moreover, although reserve requirements of CPO-NT are lower than those of TPO in general, our experiments show that: \textit{i)} reserve requirements of CPO-WR-7 are 5.28\% higher than those of TPO in quarter 2, indicating that the rule-of-thumb $\hat{\boldsymbol{r}}$ could be still insufficient in certain cases; and \textit{ii)} CPO-WR-NT leverages lower $\hat{\boldsymbol{w}}^{\diamond}$ to induce higher base-generation $P_{it}^{\star}$ in the day-ahead UC, allowing system operators to possess enough capacity and avoid load shedding when faced with $\tilde{\boldsymbol{w}}$ in the RD. Overall, this is how CPO-NT leverages the cost-oriented prediction to improve the UC economics. Last but not least, two more observations are implied by Table~\ref{tab04}:
\begin{itemize}[noitemsep, topsep=0pt, parsep=0pt, partopsep=0pt]
\item
The negative EI indicates that bi-level learning may worsen the UC economics. {Indeed, the tailors possess a common issue in wide ML applications\textemdash they cannot guarantee that learning will always bring in extra benefits. Nevertheless, this issue could be relieved with sufficient training samples and well-tuned hyper-parameters;}
\item
Although Table~\ref{tab04} indicates that the EI monotonically increases with \textit{NT}, Fig.~\ref{fig06} shows that this relationship may not hold at a 1-sample grain. Two possible reasons are responsible for this fluctuation: \textit{i)} a single tailor is used throughout one week, regardless of weekdays and weekends; and \textit{ii)} the number of training samples is limited. To relieve the fluctuation, the first approach is to train different tailors for weekdays and weekends, and the second approach is to increase \textit{NT}\textemdash the best performance is generally achieved when \textit{NT} reaches about ten times the number of parameters of the tailor. However, the two approaches would significantly increase the training time. Thus, operators should weigh up the pros and cons according to the computation budget.
\end{itemize}

Table~\ref{tab05} indicates that tailored predictions have worse MAE and RMSE. This is because the tailors evaluate the prediction quality via the actual operating cost $c^{\text{act}}$ instead of traditional accuracy metrics. Moreover, it is interesting to observe that $\hat{\boldsymbol{w}}^{\diamond}$ outperforms $\hat{\boldsymbol{w}}$ in MAPE. This is because $\hat{\boldsymbol{w}}^{\diamond}$ generally under-predicts RES, hence the MAPE of $\hat{\boldsymbol{w}}^{\diamond}$ is smaller than that of $\hat{\boldsymbol{w}}$ in the cases of scarce RES. However, these cases have a limited impact on both MAE and RMSE, because they are relatively insensitive to small values. To this end, Table~\ref{tab05} implies a fact\textemdash assessing prediction quality with a single accuracy metric may be myopic as each metric owns a specific property. \textit{To this end, this paper suggests evaluating the prediction quality with the ultimate operation goal.}

In summary, by leveraging the bi-level MIP structure to learn the UC-RD operation information, the RES-and-reserve tailor can customize raw $\{\hat{\boldsymbol{w}}, \hat{\boldsymbol{r}}\}$ to cost-oriented $\{\hat{\boldsymbol{w}}^{\diamond}, \hat{\boldsymbol{r}}^{\diamond}\}$ for improving UC economics. Although the tailored predictions are slightly worse in certain accuracy metrics, they enable a more economic operation for systems. Thus, by assessing the prediction performance via the UC economics, it is direct to conclude that tailored $\{\hat{\boldsymbol{w}}^{\diamond}, \hat{\boldsymbol{r}}^{\diamond}\}$ significantly outperforms raw $\{\hat{\boldsymbol{w}}, \hat{\boldsymbol{r}}\}$. The reduction in prediction accuracy, as listed in Table~\ref{tab05}, shall not be considered a compromise, because the ultimate goal of the prediction is to improve the optimization performance instead of simply improving the prediction accuracy. More importantly, the tailors are deployed as add-ons to the operators' accuracy-oriented predictors and thus are compatible with their current practice.

\subsection{CPO vs CPO-$\phi$}

\begin{figure}[tb]
	\centering
		\includegraphics[width=0.95\columnwidth]{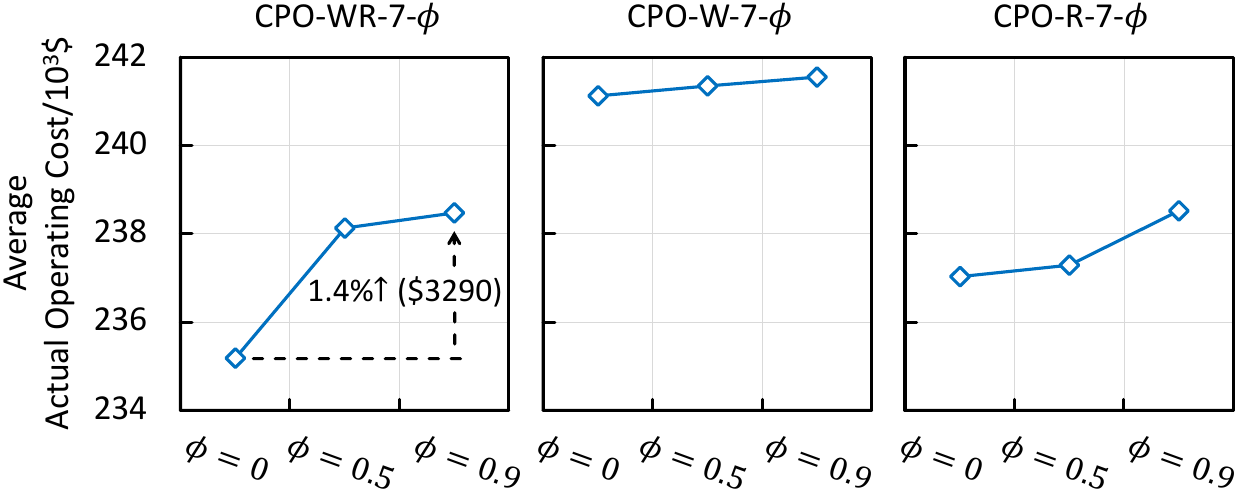}
		\vspace{-2mm}
	\caption{Out-of-sample actual operating costs of different CPO-NT-$\phi$.}\label{fig08}
	\vspace{-2mm}
\end{figure}

The binary requirement constraints (i.e., \eqref{DUC:16} and \eqref{DED:14}) are crucial in the operation models. To illustrate the impact of these constraints, this subsection compares CPO-NT-$\phi$ with different settings of the relaxing degree $\phi$.

\subsubsection{Results}
\begin{table}[tb]
	\caption{Relationship Between Relaxing Degree and ERM Objective Value (In-Sample Training Results of Quarter-4 Week)}\label{tab06}
	\centering
	\footnotesize
\begin{tabularx}{\columnwidth}{@{\extracolsep{\fill}}lccccc}
\toprule
                                       &\multicolumn{3}{c}{Objective Value of ERM \eqref{ERM}/$\text{10}^{\text{5}}$\$}\\
\cmidrule{2-4}
                                       &CPO-WR-7-$\phi$ &CPO-W-7-$\phi$   &CPO-R-7-$\phi$        \\
\midrule
$\phi=\text{0}$ (0\% Relaxation)       &3.16            &3.56             &3.15 \\
$\phi=\text{0.9}$ (90\% Relaxation)    &2.15            &2.16             &2.15 \\
\bottomrule
\end{tabularx}
\vspace{-2mm}
\end{table}

\begin{table}[tb]
	\caption{Relationship Between Relaxing Degree and Training Time (In-Sample Training Results of Quarter-4 Week)}\label{tab07}
	\centering
	\footnotesize
\begin{tabularx}{\columnwidth}{@{\extracolsep{\fill}}lccccc}
\toprule
                                       &\multicolumn{3}{c}{Training Time/s}\\
\cmidrule{2-4}
                                       &CPO-WR-7-$\phi$ &CPO-W-7-$\phi$   &CPO-R-7-$\phi$        \\
\midrule
$\phi=\text{0}$ (0\% Relaxation)       &70.81            &200.52          &81.20 \\
$\phi=\text{0.9}$ (90\% Relaxation)    &35.11            &37.63           &25.63 \\
\bottomrule
\end{tabularx}
\vspace{-3mm}
\end{table}

Fig.~\ref{fig08} compares the average actual operating costs (over the four weeks in Table~\ref{tab02}) of $\phi=\text{0}$, $\phi=\text{0.5}$, and $\phi=\text{0.9}$. It is observed that the actual operating cost increases with $\phi$, especially when the reserve requirement tailor is involved (i.e., CPO-WR-7-$\phi$ and CPO-R-7-$\phi$ in Fig.~\ref{fig08}). The largest cost increment reaches 1.4\% between CPO-WR-7-0 and CPO-WR-7-0.9, representing a daily cost increment of \$3,290. These observations indicate that tightening the binary requirement in the in-sample operation models (i.e., \eqref{ERM:4} and \eqref{ERM:5}) can enhance the out-of-sample performance of tailors.

\subsubsection{Discussions}
The binary requirement is pivotal for ensuring the clear physical meaning of decision variables within the operation models. Considering a relaxed UC model with $\phi=\text{0.9}$, to achieve a lower objective value \eqref{DUC:1}, constraint \eqref{DUC:9} allows on-off status variable $I_{it}$ and NR commitment variable $O_{it}$ to be non-zero simultaneously. That is, it allows unit $i$ to provide both SR and NR for hour $t$, which violates the physical operation rules.

Table~\ref{tab06} further reports the ERM objective value \eqref{ERM:1} (in-sample actual operating cost) of the quarter-4 training, showing that the costs of $\phi=\text{0.9}$ are significantly lower than those of $\phi=\text{0}$. However, as explained above, when $\phi=\text{0.9}$, many of the in-sample operation decisions could be infeasible in practice. As a result, such compromised in-sample operation decisions will adversely affect the out-of-sample performance of tailors. Moreover, our experimental results report that such adverse effects will be further amplified when the hyper-parameters, $\lambda^{\text{w}}$ and $\lambda^{\text{r}}$, take smaller values. For example, when $\lambda^{\text{w}}$ and $\lambda^{\text{r}}$ decrease from $\text{10}^{\text{5}}$ to $\text{10}^{\text{2}}$, the difference in actual operating costs between CPO-WR-7-0 and CPO-WR-7-0.9 increases from 1.4\% to 8.1\%.

Finally, Table~\ref{tab07} indicates that relaxing the binary requirement may accelerate the training process, because solving \textbf{MP} \eqref{MP} becomes faster. Therefore, from the perspective of operators, even though rigorously incorporating binary variables in the training is beneficial in reducing the actual operating cost, the associated long training time will necessitate a trade-off between the out-of-sample performance and training time, especially when the computation resources are limited.

\subsection{CPO vs NNPO}
The underlying idea of NNPO-NT \cite{yiwang1} is to build a differentiable function to approximate the relationship between the statistical error of a \textit{single} prediction factor (e.g., RMSE of RES predictions) and the actual operating cost, which is then used as the loss function to train NNs. 

To facilitate the comparison, the settings below are adopted:
\begin{itemize}[noitemsep, topsep=0pt, parsep=0pt, partopsep=0pt]
\item
The loss function of day $D$ is approximated by $\textit{NT}$ scenarios in day $D-1$;

\item
NNPO-NT uses tailored RES prediction and keeps raw reserve requirements. This is because NNPO-NT cannot be directly applied to tailor multiple prediction factors if these factors affect each other mutually. Therefore, for the sake of fairness, NNPO-NT is compared to CPO-W-NT;

\item
Dec 15-18, 2020, the last four days in Table~\ref{tab02}, are used for comparison. Hence, up to 349 training samples (from Jan-1 to Dec-14, 2020) are available for NNPO-NT;

\item
Two NNPO-NTs\textemdash NNPO-7 and NNPO-349\textemdash are trained for comparison, in which NNPO-7 uses the same training samples as CPO-W-7;

\item
Note that randomness exists in the NN training due to certain heuristic steps. Thus, NNPO-NT is evaluated statistically. Taking NNPO-7 as an example, the evaluation process involves three steps: \textit{i)} training NNPO-7 ten times, which yields ten NNPO-7s; \textit{ii)} conducting the out-of-sample test for each of the ten NNPO-7s; and \textit{iii)} calculating the mean, standard deviation, and confidence interval (CI) of their out-of-sample results, thus providing a statistical evaluation of NNPO-7.
\end{itemize}

\subsubsection{Results}
\begin{figure}[tb]
	\centering
		\includegraphics[width=\columnwidth]{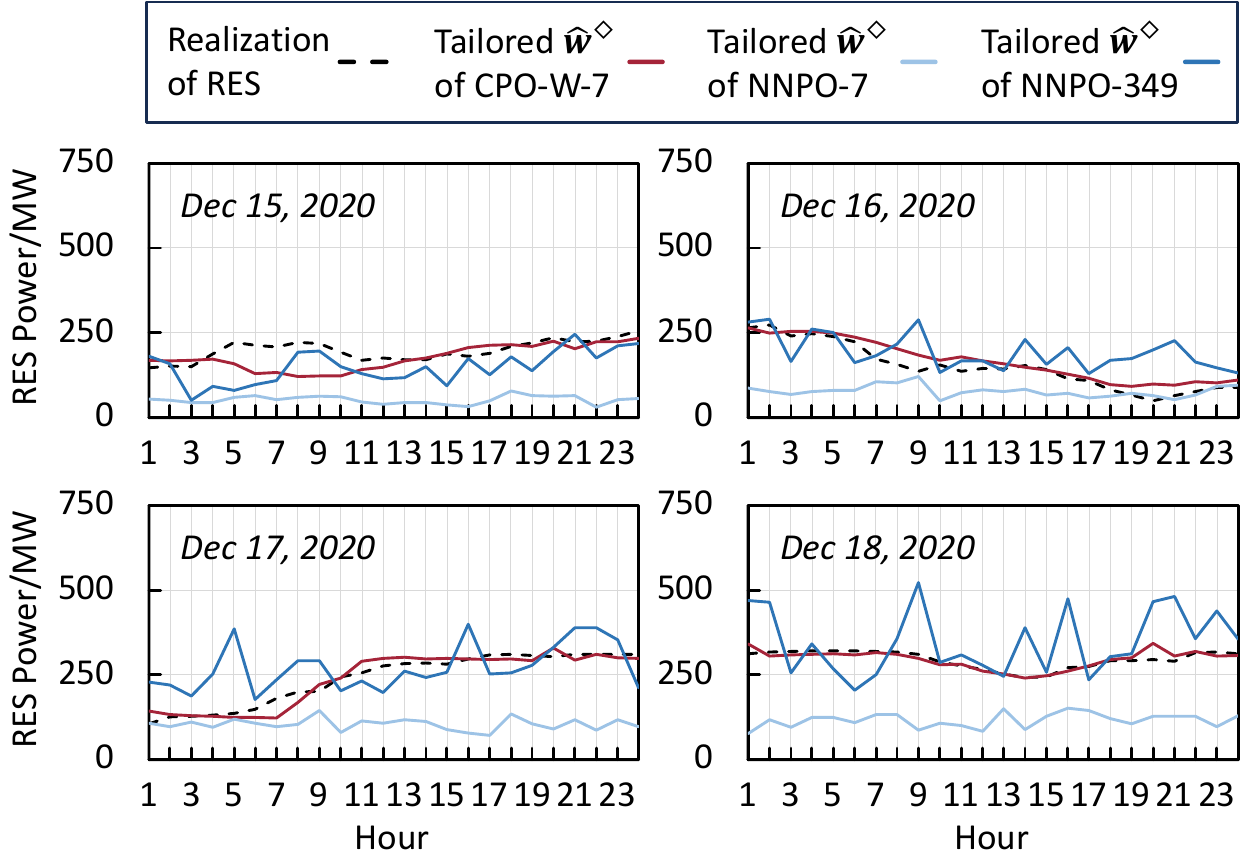}
		\vspace{-5mm}
	\caption{Tailored RES predictions of CPO-W-7 and two NNPO-NTs.}\label{fig09}
\end{figure}

Fig.~\ref{fig09} compares the tailored RES predictions, showing that: \textit{i)} NNPO-7 tends to under-predict, and its curves are fairly flat; \textit{ii)} compared to NNPO-7, the curves of NNPO-349 are higher and more volatile; and \textit{iii)} the curves of CPO-W-7 are less volatile and closer to the actual realization.

Table~\ref{tab08} lists the actual UC cost and balancing cost, showing that CPO-W-7 falls between NNPO-7 and NNPO-349 regarding the mean actual UC cost. This phenomenon can be explained by Fig.~\ref{fig09}\textemdash NNPO-7 tends to tailor the raw $\hat{\boldsymbol{w}}$ to be a flat and under-predicted $\hat{\boldsymbol{w}}^{\diamond}$, while NNPO-349 generally produces a $\hat{\boldsymbol{w}}^{\diamond}$ involving over-predicted spikes. Consequently, NNPO-7/NNPO-349 turns on more/fewer units in the day-ahead stage, causing higher/lower actual UC costs. Although the differences in the actual UC cost seem minor, they affect the balancing cost considerably. For example, on Dec 18, 2020, the mean actual UC cost of NNPO-7 is only $\$\text{0.4}\times\text{10}^{\text{3}}$ higher than that of  CPO-W-7, but the difference it causes in the mean balancing cost is as high as $\$\text{63.9}\times\text{10}^{\text{3}}$.

Table~\ref{tab08} also indicates that balancing costs of NNPO-NT have noticeable standard deviations, NNPO-7 in particular. This is because the gradient-based NN training involves certain heuristic steps and can only provide a locally optimal tailor, i.e., it cannot guarantee a consistent training result. In contrast, CPO-W-7 does not exhibit such noticeable standard deviations, as its training relies on mathematical programming that can consistently approach a globally optimal predictor.

\begin{table}[tb]
	\caption{Out-of-Sample Costs of CPO-W-7 and Two NNPO-NTs}\label{tab08}
	\centering
	\footnotesize
\begin{tabularx}{\columnwidth}{@{\extracolsep{\fill}} llccccc}

\toprule
                       &          &\multicolumn{2}{c}{Actual UC Cost/10$^{\text{3}}$\$}     &\multicolumn{2}{c}{Balancing Cost/10$^{\text{3}}$\$} \\
\cmidrule{3-4} \cmidrule{5-6}
                       &          &Mean & Std. Dev.                     &Mean  & Std. Dev. \\
\midrule
\multirow{3}{*}{Dec 15}& CPO-W-7  &3.1  &0.0                            &223.1 &0.0         \\
                       & NNPO-7   &4.2  &0.4                            &264.9 &5.8         \\
                       & NNPO-349 &3.8  &0.2                            &231.6 &2.4         \\
\midrule
\multirow{3}{*}{Dec 16}& CPO-W-7  &3.6  &0.0                            &250.9 &0.0         \\
                       & NNPO-7   &4.1  &0.3                            &284.3 &8.0         \\
                       & NNPO-349 &3.5  &0.3                            &256.1 &3.2         \\
\midrule
\multirow{3}{*}{Dec 17}& CPO-W-7  &2.9  &0.0                            &190.8 &0.0            \\
                       & NNPO-7   &3.5  &0.2                            &228.3 &4.5            \\
                       & NNPO-349 &2.8  &0.2                            &198.0 &3.4            \\
\midrule
\multirow{3}{*}{Dec 18}& CPO-W-7  &2.9  &0.0                            &146.6 &0.0            \\
                       & NNPO-7   &3.3  &0.0                            &210.5 &8.2            \\
                       & NNPO-349 &2.7  &0.2                            &161.1 &5.2            \\
\bottomrule

\end{tabularx}
\end{table}

\begin{figure}[tb]
	\centering
		\includegraphics[width=\columnwidth]{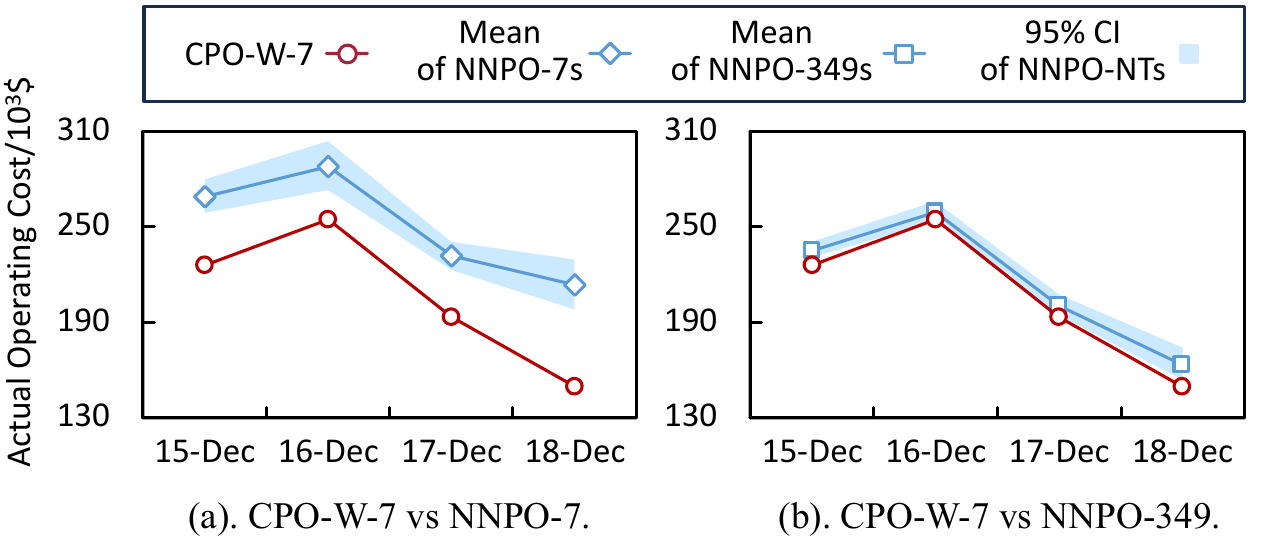}
		\vspace{-6mm}
	\caption{Out-of-sample operating costs of CPO-W-7 and two NNPO-NTs.}\label{fig10}
\end{figure}

Fig.~\ref{fig10}(a) shows that using the same 7 training samples, CPO-W-7 significantly outperforms NNPO-7 on the actual operating cost. This finding indicates that CPO-W-7 is more effective than NNPO-7 when training samples are limited\textemdash a benefit attributed to the bi-level structure of the ERM problem \eqref{ERM}. Furthermore, Fig.~\ref{fig10}(b) shows that compared to NNPO-7, NNPO-349 has a lower mean and a narrower CI for the actual operating cost. However, CPO-W-7 still economically outperforms NNPO-349.

Overall, the above results have clearly demonstrated that, even though the number of training samples is limited, CPO can still learn to properly tailor the raw $\hat{\boldsymbol{w}}$. The tailored $\hat{\boldsymbol{w}}^{\diamond}$ leads to a moderate actual UC cost and a low balancing cost, ultimately resulting in a low actual operating cost.

\subsubsection{Discussions}
To facilitate a more in-depth comparison between CPO and NNPO, their differences and the underlying reasons are further discussed.

\textit{\textbf{Performances with Few Training Samples:}}
Our experiments report that NNPO-7 performs poorly in both in-sample and out-of-sample losses\textemdash a clear sign of under-fitting. This is mainly due to the existence of ``bad'' training samples. With only 7 training samples, even one or two poor-quality samples could noticeably mislead the training and worsen the generalizability of NNPO-7. Our experiments also report that although NNPO-7 can be improved by removing the ``bad'' samples, it remains economically worse than CPO-W-7.

In comparison, CPO-W-7 exhibits robustness against the ``bad'' training samples. This advantage primarily arises from the bi-level structure of \eqref{ERM}. Within this structure, the non-anticipativity between the prediction, UC, and RD stages can be ensured, and the UC and RD processes can be modeled uncompromisingly. As a result, such an ERM model can accurately express the predict-then-optimize practice\textemdash sequentially conduct predictions, UC optimization \eqref{UC}, and RD optimization \eqref{ED} with the information at hand (i.e., non-anticipativity). Consequently, while tailors trained via the bi-level ERM might not show a considerably low in-sample loss, they can still achieve good out-of-sample performance.

\textit{\textbf{Tuning Burdens of Hyper-Parameters:}}
The high performance of NNPO can be achieved by tuning hyper-parameters finely. However, tuning a vast array of hyper-parameters on a daily basis is challenging. Indeed, most commercial NN-based predictors, once delicately tuned, rarely undergo frequent re-tuning processes due to significant computational burdens. As compared, CPO only requires tuning the two hyper-parameters, $\lambda^{\text{w}}$ and $\lambda^{\text{r}}$. Thus, from the perspective of hyper-parameter tuning, CPO is more user-friendly.

\textit{\textbf{Interpretability of Tailored Predictions:}}
It is noteworthy that the NN-based tailor may yield \textit{negative} RES predictions. This issue, together with NNPO-349's spikes observed in Fig.~\ref{fig09}, could raise operators' concerns: \textit{i)} dramatically tailoring the high-quality raw predictions lacks interpretability; and \textit{ii)} the spikes would trigger the needs of extra flexibility resources. Indeed, these issues stem from the fact that the NN training ignores critical physical information (e.g., ramping and generation limits) within the UC and RD models. That is, the training exclusively aims to minimize the training loss, regardless of physical principles, e.g., RES predictions shall be non-negative.

Regarding this aspect, CPO possesses two valuable properties: \textit{i)} the linear form \eqref{form} offers strong interpretability; and \textit{ii)} the bi-level ERM \eqref{ERM} effectively integrates the UC-RD process into the training. Within this ERM, irrational predictions (e.g., spikes and negative values) will cause extreme costs and even infeasibility in the lower-level UC and/or RD, thus being excluded. As a result, tailors trained by the physically-informed ERM barely produce irrational predictions.

\textit{\textbf{Form of Tailors' Model:}}
NNPO possesses the ability to capture complex non-linear relationships. As compared, the presented tailors \eqref{form} can only scale raw predictions linearly. Nevertheless, for daily UC practice, tailors with simpler forms may be more appropriate. This is because operators have already utilized sophisticated models like NNs to produce high-quality raw predictions \cite{Intro_RESPRE}; thus, only moderate fine-tuning instead of dramatic modification would be needed. Indeed, our extensive experiments show that NNPO performs better with linear activation than with ReLU, sigmoid, softmax, and tanh activations. This observation implies that straightforward tailor models could be more effective when high-quality raw predictions are used as inputs.

\subsection{CPO vs TSP and TRO Models with MIP Recourse}
CPO-WR-NT, TSP-NS, and TRO-NB are further compared using a modified IEEE 118-bus system (with 8,600MW non-RES capacity and 1,000MW RES capacity). For the sake of comparison, the following settings on TSP-NS \eqref{SPUC} and TRO-NB \eqref{ROUC} models with MIP recourse are adopted:
\begin{itemize}[noitemsep, topsep=0pt, parsep=0pt, partopsep=0pt]
\item
The reserve requirement constraints \eqref{DUC:19} are not explicitly included in the two-stage models \eqref{SPUC}-\eqref{ROUC}, because the reserve schedules have been implicated in their second-stage RD. To enable the first-stage UC to deploy enough reserve, not only RES but also load are regarded as uncertainty factors. The day-ahead raw point predictions (denoted as $\hat{\boldsymbol{u}}$) and their corresponding confidence intervals are provided by the Belgian system \cite{Intro_RESPRE};
\item
TSP-NS is modeled as a scenario-based model \eqref{SPUC}. First, 3,000 scenarios are generated via Latin hypercube sampling within the 90\% confidence intervals. And then, these 3,000 scenarios are reduced to \textit{NS} scenarios via a scenario-tree method \cite{tree}. For scenario $h$, $\tilde{\boldsymbol{u}}^{\text{scan}}_{h}$ denotes its realization, and $p_{h}$ denotes its probability. Model \eqref{SPUC} is directly solved by Gurobi;
\begin{subequations}\label{SPUC}
\begin{align}
\min_{\boldsymbol{x}, \boldsymbol{y}, \boldsymbol{z}}\,
        &  \boldsymbol{b}^{\top}\boldsymbol{x} + \boldsymbol{c}^{\top}\boldsymbol{y} + p_{h} \textstyle{\sum\limits_{h = 1}^{\text{\textit{NS}}} \boldsymbol{d}^{\top}\boldsymbol{z}_{h}}                                       \label{SPUC:1}\\
s.\,t.\,& \boldsymbol{x}, \boldsymbol{y} \in
                     \mathcal{X}(\hat{\boldsymbol{u}});                                        \label{SPUC:2}\\
        & \boldsymbol{z}_{h} \in
                     \mathcal{Z}(\boldsymbol{x},\boldsymbol{y},
                     \tilde{\boldsymbol{u}}^{\text{scen}}_{h}),\, h = 1,..., \text{\textit{NS}};   \label{SPUC:3}
\end{align}
\end{subequations}

\item
TRO-NB, as in \eqref{ROUC}, uses an adjustable budget parameter $\Gamma \in [\text{0}, \text{48}]$ to control the solution robustness, in which $\Gamma = \text{\textit{NB}}$. The box uncertainty set $\mathcal{U}$ is built using the 90\% confidence intervals, where the worst-case realization of the uncertainty is denoted as $\tilde{\boldsymbol{u}}^{\text{wcr}}$. A nested C\&CG method \cite{RO1} is used to solve \eqref{ROUC}.
\begin{subequations} \label{ROUC}
\begin{align}
\min_{\boldsymbol{x},\boldsymbol{y}}
              & \,  \boldsymbol{b}^{\top}\boldsymbol{x} + \boldsymbol{c}^{\top}\boldsymbol{y} + \max_{\tilde{\boldsymbol{u}}^{\text{wcr}} \in \mathcal{U}} \min_{\boldsymbol{z}}\boldsymbol{d}^{\top}\boldsymbol{z}   \label{ROUC:1}\\
s.\,t.\,      &  \boldsymbol{x},\boldsymbol{y} \in \mathcal{X}(\hat{\boldsymbol{u}});\,       \label{ROUC:2}\\
              &  \boldsymbol{z} \in \mathcal{Z}(\boldsymbol{x},\boldsymbol{y},
                                         \tilde{\boldsymbol{u}}^{\text{wcr}},
                                         \Gamma);                                             \label{ROUC:3}
\end{align}
\end{subequations}
\end{itemize}

\begin{figure}[tb]
	\centering
		\includegraphics[width=\columnwidth]{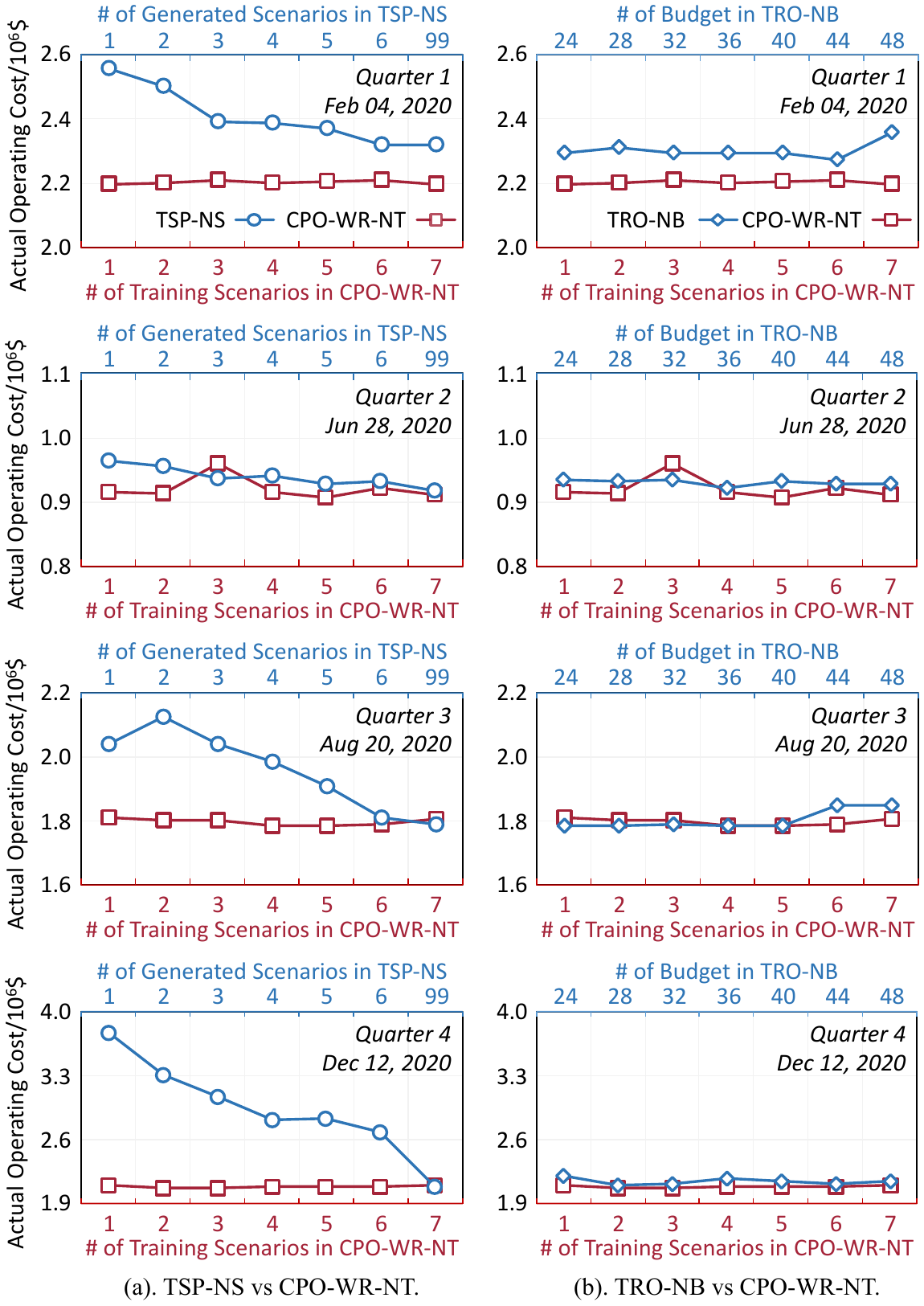}
		\vspace{-5mm}
	\caption{Comparison of CPO-WR-NT (red) and TSP-NS/TRO-NB (blue).}\label{fig11}
\end{figure}

\subsubsection{Results of CPO vs TSP-NS}
Fig.~\ref{fig11}(a) sketches the UC economics of seven CPO-WR-NTs and seven TSP-NSs, leading to the following three observations:
\begin{itemize}[noitemsep, topsep=0pt, parsep=0pt, partopsep=0pt]
\item
TSP-NS requires many scenarios to perform well, leaving a wide cost gap between TSP-1 and TSP-99. Relatively, CPO-WR-NT needs fewer samples, with a narrower cost gap between CPO-WR-1 and CPO-WR-7. This point means that CPO can leverage training samples more effectively, ensuring its performance in the case of few training samples;

\item
Except for quarter 1, the actual operating cost of TSP-NS and CPO-WR-NT eventually intersects, implying that the performances of TSP-NS and CPO-WR-NT are mutually reachable with suitable \textit{NS}/\textit{NT};

\item
In quarter 1, even with \textit{NS} = 99, TSP-NS still performs poorly. This is because the confidence intervals cannot exactly cover the underlying distributions of $\tilde{\boldsymbol{u}}$, and TSP-NS would encounter significant slack penalties.
\end{itemize}

Table~\ref{tab09} breaks down the actual operating cost of the quarter-2 day. It shows that the three TSP-NSs bring lower startup and no-load costs. This is because the reserve requirements \eqref{DUC:19} are implicated in \eqref{SPUC} rather than being explicit as in \eqref{PUC}. Thus, the three TSP-NSs turn on fewer units and schedule less reserve in the UC stage, causing relatively more extensive usage of quick-start units in the RD stage. Indeed, our results show that the average hourly reserve schedules (SR plus NR) of TSP-1 and TSP-99 are, respectively, 89.9MW and 207.2MW. In comparison, these average schedules are 1,265.8MW and 1,286.8MW in CPO-WR-1 and CPO-WR-7. As quick-start units have limited ramping ranges, they are operated away from the most efficient points in RD, causing higher generation costs in TSP-NS.
\begin{table}[tb]
	\caption{Breakdown of Actual Operating Cost on 28 Jun 2020}\label{tab09}
	\centering
	\footnotesize
\begin{tabularx}{\columnwidth}{@{\extracolsep{\fill}} lcccccc}
\toprule
        &\multicolumn{2}{c}{Actual UC Cost/10$^{\text{3}}$\$}& \multicolumn{3}{c}{Balancing Cost/10$^{\text{3}}$\$}\\
       \cmidrule{2-3} \cmidrule{4-6}
        &\multirow{2}{*}{Startup}&\multirow{2}{*}{No-Load}&Startup      &\multirow{2}{*}{Gen.}&\multirow{2}{*}{Slack}\\
        &                        &                        &plus No-Load &                     &                \\
\midrule
CPO-WR-7&6.7                     & 9.0                    & 0.0         & 894.9               & 00.0     \\
TSP-1   &6.5                     & 8.1                    & 1.1         & 934.3               & 13.7     \\
TSP-4   &5.1                     & 7.7                    & 1.0         & 926.4               & 00.0     \\
TSP-99  &4.8                     & 7.8                    & 0.7         & 904.7               & 00.0     \\
\bottomrule
\end{tabularx}
\end{table}

Indeed, the difference in reserve schedules implies that TSP-NS and CPO-WR-NT achieve similar economics in different ways. Specifically, TSP-NS schedules fewer reserves and uses more RES curtailments in place of SR in the RD. This may work well in case of ample RES but becomes economically risky when the actual RES availability is low, as expensive slack penalties will be triggered (e.g., due to insufficient capacities, TSP-1 is penalized as in Table~\ref{tab09}). This may also cause a lower RES utilization.

\subsubsection{Results of CPO vs TRO-NB}
Fig.~\ref{fig11}(b) sketches the economics of CPO-WR-NT and TRO-NB, showing that both methods are relatively stable\textemdash the gap between the best and worst costs is small. Moreover, multiple intersections exist in the curves, implying that TRO-NB and CPO-WR-NT can reach similar performance by tuning parameters properly.

Table~\ref{tab10} further compares the costs on the quarter-3 day. It shows that TRO-NB has a higher startup cost but a lower no-load cost in the UC stage. The reason is that TRO-NB schedules SR for those hours of potential large volatility, so that multiple units in the UC stage are turned on before these hours and turned off afterward. On the other hand, CPO-WR-7 constantly keeps more units online to meet the reserve requirements of individual hours, with higher no-load costs. In addition, the penalty of TRO-48 means that even a large $\Gamma$ could lead to an insufficient reserve schedule. Indeed, the hourly average reserve schedules of TRO-48 and CPO-WR-7 are 278.4MW and 1,214.7MW, implying that CPO-WR-7 is more conservative. With this, it can be pointed out that CPO-WR-NT achieves economics similar to the state-of-the-art methods in a more reliable way\textemdash scheduling more reserves to avoid potential slack penalties.
\begin{table}[tb]
	\caption{Breakdown of Actual Operating Cost on 20 Aug 2020}\label{tab10}
	\centering
	\footnotesize
\begin{tabularx}{\columnwidth}{@{\extracolsep{\fill}} lcccccc}
\toprule
       &\multicolumn{2}{c}{Actual UC Cost/$\text{10}^{\text{3}}$\$}
       &\multicolumn{3}{c}{Balancing Cost/$\text{10}^{\text{3}}$\$}\\
       \cmidrule{2-3} \cmidrule{4-6}
        &\multirow{2}{*}{Startup}&\multirow{2}{*}{No-Load}&Startup      &\multirow{2}{*}{Gen.}&\multirow{2}{*}{Slack}\\
        &                        &                        &plus No-Load &                     &                \\
\midrule
CPO-WR-7  &5.9        &12.6       &0.2                   &1789.7         &0\\
TRO-24 &6.5        &11.3       &1.2                   &1764.9         &0\\
TRO-36 &6.7        &11.4       &1.2                   &1764.3         &0\\
TRO-48 &6.5        &11.3       &1.2                   &1764.9         &63.5\\
\bottomrule
\end{tabularx}
\end{table}

\subsubsection{Discussions}
The above results have shown that CPO is comparable to TSP and TRO regarding UC economics. Although a clear-cut economic advantage of CPO has yet to be observed, it does offer two promising practical advantages.

\textit{\textbf{Computational Affordability to System Operators:}}
Computation time of the proposed CPO framework includes the training time of tailors and the solving time of UC tasks. As the tailors are trained offline and rather infrequently, their computational time is not considered a bottleneck for practical applications. In addition, the training time is significantly impacted by the efficiency of the codes, such as Python's TensorFlow library used in case studies. Therefore, this subsection focuses on the solving time of UC models.

The day-ahead UC to determine optimal operation plans has been regarded as one of the most computationally challenging tasks, even with the deterministic UC formulation \eqref{UC}. Table~\ref{tab11} compares the solving time of UC of TPO \eqref{UC}, CPO-WR-7 \eqref{PUC}, TSP-99 \eqref{SPUC}, and TRO-48 \eqref{ROUC}. TSP-99 and TRO-48 are much more computationally taxing than TPO. Although advanced accelerating techniques, e.g., parallel computation and decomposition methods, could accelerate the solving time of TSP and TRO, they still lag behind the deterministic UC \cite{UCspeed}. This is one of the main reasons that most operators still use the deterministic UC \eqref{UC}. As compared, the prescriptive UC \eqref{PUC} of CPO keeps a computation burden similar to \eqref{UC} of TPO, thus being more affordable to operators.
\begin{table}[tb]
	\caption{Solving Time of UC Models}\label{tab11}
	\centering
	\footnotesize
\begin{tabularx}{\columnwidth}{@{\extracolsep{\fill}}lccccc}
\toprule
&TPO      &CPO-WR-7 &TSP-99   &TRO-48        \\
\midrule
&4.21s    &4.29s    &1,209.32s &3,082.86s \\
\bottomrule
\end{tabularx}
\end{table}

\textit{\textbf{Compatibility with System Operators' Practice:}}
In practice, system operators aim to identify the economically optimal UC decisions given the predictions. However, UC decisions of TSP/TRO may not fully follow the least-cost principle in the short term due to the treatment of uncertainties. Although better UC economics may be reached in the long term, a short-term financial deficit would be inevitable \cite{vladimir}.

While pursuing to maximize system economics, operators must have a stably interpretable operation plan (e.g., which units to dispatch), so that both operation schedules and post-operation analyses can be rationally understood. Regarding this, TSP and TRO could expose intrinsic obstacles in practical applications. Indeed, for TSP, a stable schedule result requires sufficient scenarios, which, however, will further magnify the computational challenge; while for TRO, to assure a tractable reformulation, the uncertainty set is usually designed to be heavily parameter-dependent. These together cause volatility and poor interpretability of schedule results in TSP and TRO.

Specifically, the process of constructing scenario sets for TSP inevitably introduces randomness. However, different scenario sets, even with tiny differences, may lead to noticeable changes in commitment results, i.e., certain units are ON/OFF from one set to another. These issues adversely affect the interpretability of the schedule plan. Analogously, the setting of parameter $\Gamma$ suffers from the same issue. This is because the impacts of different $\Gamma$ on the schedule plan are difficult to interpret and quantify. Thus, the operators have to face the same issue on a daily basis\textemdash what is the best scenario set or budget parameter for the next operation day?

Last but not least, the difference between the scheduled reserves of CPO and TSP/TRO implies that the second stage of TSP/TRO may not be able to implicitly secure enough reserves in the UC plan, indicated by the non-zero slacking penalties in Tables \ref{tab09} and \ref{tab10}. In other words, TSP and TRO may not be as robust as expected against uncertainties, implying that the explicit expression \eqref{DUC:19} of reserve requirements would remain necessary for the UC model. Regarding this, CPO leverages the raw predictions of reserve requirements to achieve good robustness.

In summary, the prescriptive UC \eqref{PUC} of CPO has the same structure as the model \eqref{UC} of TPO, thus following the least-cost principle, remaining strong interpretability, and being computationally efficient. Moreover, Fig.~\ref{fig11} has shown that the UC economics of TSP/TRO and CPO are mutually reachable. Thus, by selecting training samples properly, it is promising to harness the compatibility of CPO in bridging the gap between TSP/TRO and current practice.

\section{Conclusion}\label{sec05}
To improve the existing predict-then-optimize practice and reduce the actual operating cost, this paper presents a CPO framework based on the bi-level MIP formulation. The idea is to train an add-on tailor that can customize raw predictions of RES and reserve requirements into cost-oriented predictions, thereby reducing the actual operating cost. Case studies compare the presented framework to five other frameworks in the literature, leading to the following conclusions:
\begin{itemize}[noitemsep, topsep=0pt, parsep=0pt, partopsep=0pt]
\item
Compared to the existing predict-then-optimize practice, results over four selective weeks from different seasons demonstrate that the presented framework can reduce the average actual operating cost by up to 2.54\%;

\item
Compared to the bi-level relaxed MIP method, the presented framework with an uncompromising binary requirement can achieve up to an 8.1\% reduction in average actual operating costs;

\item
Compared to the NN-based CPO framework, the presented bi-level MIP-based CPO framework exhibits two advantages: \textit{i)} it can reduce the actual operating cost even when the training samples are few; and \textit{ii)} it has a relatively consistent and stable out-of-sample performance;

\item
The actual operating cost induced by the presented framework, TSP, and TRO could be mutually reachable. However, since the RES-and-reserve tailor acts as an add-on to the deterministic UC, the presented framework offers higher implementation values to the predict-then-optimize practice.
\end{itemize}

Future work could further consider other UC-related issues, including: \textit{i)} how can the impacts of system topology changes (e.g., transmission contingency and transmission line switching) on tailor training be managed? \textit{ii)} what strategies can be employed to effectively handle data drifts (e.g., seasonal and weekly load pattern changes)? and \textit{iii)} how can the computational efficiency of the proposed solution approach be further enhanced to accommodate larger-scale systems?

\appendix
\subsection{The Deterministic UC Model with Predictions}\label{app_a}
The UC objective \eqref{DUC:1} is to minimize the total system cost, including startup, no-load, and generation costs. The generator constraints include generation limits \eqref{DUC:2}-\eqref{DUC:3}, segment-based generation representation \eqref{DUC:4}-\eqref{DUC:5}, SR \eqref{DUC:6} and NR \eqref{DUC:7}-\eqref{DUC:9} capacity limits, startup-shutdown-commitment status logic \eqref{DUC:10}, minimum on and off requirements \eqref{DUC:11}-\eqref{DUC:12}, ramping limits \eqref{DUC:13}-\eqref{DUC:14}, and RES power limits \eqref{DUC:15} \cite{zl1}. Constraint \eqref{DUC:16} is the binary requirement. System constraints include power balance \eqref{DUC:17}, transmission limits \eqref{DUC:18} based on DC power flow-based function $\mathcal{F}(\cdot)$ \cite{Cost}, and reserve requirements \eqref{DUC:19}. Note that $\hat{L}_{q}$ denotes the load demand prediction of bus $q$.
\begin{subequations}\label{DUC}
\begin{flalign}
&\min_{\boldsymbol{x}, \boldsymbol{y}}
\textstyle{\sum\limits_{i \in \mathcal{I}} \sum\limits_{t \in \mathcal{T}}
(C_{i}^{\text{su}}U_{it} + C_{i}^{\text{nl}}I_{it})
+ \sum\limits_{i \in \mathcal{I}} \sum\limits_{t \in \mathcal{T}} \sum\limits_{k \in \mathcal{K}} C^{\text{sg}}_{ik}P_{itk}^{\text{sg}}}                \mspace{-350mu}&                                                 \label{DUC:1}        \\
&\text{where}\,\,\boldsymbol{x} = \{\boldsymbol{D}, \boldsymbol{I}, \boldsymbol{O}, \boldsymbol{U} \}, \boldsymbol{y} = \{\boldsymbol{P}, \boldsymbol{P}^{\text{sg}}, \boldsymbol{R}^{\text{sr}}, \boldsymbol{R}^{\text{nr}}, \boldsymbol{W}\}  \mspace{-350mu}&                                                        \notag\\
&\textit{Generator Constraints:}\mspace{-350mu}&                                                        \notag        \\
&\textstyle{P_{it}-R_{it}^{\text{sr}}\geq P^{\text{m}}_{i}I_{it},}
                               \mspace{-350mu}&\forall i\in \mathcal{I}, t\in \mathcal{T};              \label{DUC:2}\\
&\textstyle{P_{it}+R_{it}^{\text{sr}}\leq P^{\text{M}}_{i}I_{it},}
                               \mspace{-350mu}&\forall i\in \mathcal{I}, t\in \mathcal{T};              \label{DUC:3}\\
&\textstyle{P_{it} = \sum_{k \in \mathcal{K}} P_{itk}^{\text{sg}},}
                               \mspace{-350mu}&\forall i\in\mathcal{I}, t\in\mathcal{T};                \label{DUC:4}\\
&\textstyle{0 \leq P_{itk}^{\text{sg}} \leq \bar{P}_{ik}^{\text{sg}} I_{it},}
                               \mspace{-350mu}&\forall i\in\mathcal{I}, t\in\mathcal{T},k\in\mathcal{K};\label{DUC:5}\\
&\textstyle{0 \leq R_{it}^{\text{sr}} \leq \bar{R}_{i}^{\text{sr}}I_{it}},
                               \mspace{-350mu}&\forall i\in \mathcal{I}, t\in \mathcal{T};              \label{DUC:6}\\
&\textstyle{P^{\text{m}}_{i}O_{it}\leq R_{it}^{\text{nr}}\leq \bar{R}_{i}^{\text{nr}}O_{it},}
                               \mspace{-350mu}&\forall i\in \mathcal{I}, t\in \mathcal{T};              \label{DUC:7}\\
&\textstyle{O_{it} = 0,}       \mspace{-350mu}&\forall i\in \mathcal{I}^{\text{ns}},t\in \mathcal{T};   \label{DUC:8}\\
&\textstyle{O_{it}+I_{it}\leq 1},\mspace{-350mu}&\forall i\in \mathcal{I}, t\in \mathcal{T};            \label{DUC:9}\\
&\textstyle{U_{it}-D_{it}=I_{it}-I_{i,t-1},}
                               \mspace{-350mu}&\forall i\in \mathcal{I}, t\in \mathcal{T};              \label{DUC:10}\\
&\textstyle{\sum_{t^\prime=t-T_{i}^{\text{su}}+1}^{t} U_{i t^\prime}\leq I_{it}},
                               \mspace{-350mu}&\forall i\in\mathcal{I},t\in\mathcal{T}^{\text{su}}_{i}; \label{DUC:11}\\
&\textstyle{\sum_{t^\prime=t-T_{i}^{\text{sd}}+1}^{t} D_{i t^\prime}\leq 1 - I_{it},}
                               \mspace{-350mu}&\forall i\in\mathcal{I},t\in\mathcal{T}^{\text{sd}}_{i}; \label{DUC:12}\\
&\textstyle{P_{it}-P_{i,t-1}\leq P^{\text{M}}_{i}(1-I_{it})},
                               \mspace{-350mu}&                                                         \notag       \\
&\textstyle{\mspace{30mu} + R_{i}^{\uparrow}I_{i,t-1} + R_{i}^{\text{su}}(I_{it}- I_{i,t-1})},
                               \mspace{-350mu}&\forall i\in \mathcal{I},t\in \mathcal{T};               \label{DUC:13}\\
&\textstyle{P_{i,t-1} - P_{it}\leq P^{\text{M}}_{i} (1 - I_{i,t-1})},
                               \mspace{-350mu}&                                                         \notag        \\
&\textstyle{\mspace{30mu} + R^{\downarrow}_{i} I_{it} + R^{\text{sd}}_{i}(I_{i,t-1} - I_{it})},
                               \mspace{-350mu}&\forall i\in \mathcal{I},t\in \mathcal{T};               \label{DUC:14}\\
&\textstyle{0\leq W_{jt}\leq\hat{{W}}_{jt},}
                               \mspace{-350mu}&\forall j\in \mathcal{J},t\in \mathcal{T};               \label{DUC:15}\\
&\textstyle{D_{it},I_{it},O_{it},U_{it}\in\{0,1\},}
                               \mspace{-350mu}&\forall i\in \mathcal{I},t\in \mathcal{T};               \label{DUC:16}\\
&\textit{System Constraints:}  \mspace{-350mu}&                                                         \notag       \\
&\textstyle{\sum\nolimits_{i \in \mathcal{I}} P_{it}
         + \sum\nolimits_{j \in \mathcal{J}} W_{jt} = \sum\nolimits_{q \in \mathcal{Q}} \hat{L}_{qt},}
                               \mspace{-350mu}&\forall t\in \mathcal{T};                                \label{DUC:17}\\
&\textstyle{-B_{b}\leq \mathcal{F}_{b}(\boldsymbol{P}_{t}, \boldsymbol{W}_{t}, \hat{\boldsymbol{L}}_{t}) \leq B_{b},}
                               \mspace{-350mu}&\forall b\in \mathcal{B}, t\in \mathcal{T};              \label{DUC:18}\\
&\textstyle{\sum\limits_{i \in \mathcal{I}}R_{it}^{\text{sr}} \geq \hat{R}_{t}^{\text{sr}},\,\sum\limits_{i \in \mathcal{I}}(R_{it}^{\text{sr}} + R_{it}^{\text{nr}}) \geq \hat{R}_{t}^{\text{sr}} + \hat{R}_{t}^{\text{nr}},}
                               \mspace{-350mu}&\forall t\in \mathcal{T};                                \label{DUC:19}
\end{flalign}
\end{subequations}

\subsection{The Re-Dispatch Model with Realizations}\label{app_b}
In the RD model \eqref{DED}, those variables with superscript $\star$ are optimal variables provided by solving \eqref{DUC}. The objective \eqref{DED:1} is to minimize startup and no-load costs of quick-start units that are not committed in UC but scheduled for providing NR, generation costs of all units, and slack penalty costs. Constraints \eqref{DED:2}-\eqref{DED:4} mean that only the quick-start units not committed in UC may change statuses. All units are subjected to generation limits \eqref{DED:5}-\eqref{DED:8}, dispatch adjustment limits \eqref{DED:9}-\eqref{DED:10}, ramping limits \eqref{DED:11}-\eqref{DED:12}, RES power limits \eqref{DED:13}, and binary requirements \eqref{DED:14}. The system shall satisfy power balance \eqref{DED:15} and transmission limits \eqref{DED:16}-\eqref{DED:17} against actual load $\tilde{L}$ and RES realizations. The non-negative slack variables \eqref{DED:18}-\eqref{DED:19} ensure the RD feasibility w.r.t. the given UC solutions.
\begin{subequations}\label{DED}
\begin{flalign}
& \min_{\boldsymbol{z}}
\left\{ \begin{array}{l} 
      \sum\limits_{i \in \mathcal{I}} \sum\limits_{t \in \mathcal{T}}
      (C_{i}^{\text{su}}U_{it}^{\text{RD,qs}}+C_{i}^{\text{nl}}I_{it}^{\text{RD,qs}}) \\
     \quad+ \sum\limits_{i \in \mathcal{I}} \sum\limits_{t \in \mathcal{T}}\sum\limits_{k \in \mathcal{K}}C_{ik}^{\text{sg}}P_{itk}^{\text{RD,sg}} \\
     \quad+\sum\limits_{t \in \mathcal{T}} (C^{\text{gs}}S_{t}^{\text{gs}} + C^{\text{ls}}S_{t}^{\text{ls}})\\
     \quad+\sum\limits_{t \in \mathcal{T}}\sum\limits_{b \in \mathcal{B}}C^{\text{bs}}(S_{bt}^{+} + S_{bt}^{-})
     \end{array} \right\}
                                                \mspace{-290mu}&   \label{DED:1}\\
&\text{where}\,\boldsymbol{z}=
\left\{ \begin{array}{l}
\boldsymbol{D}^{\text{RD,qs}},
\boldsymbol{I}^{\text{RD}},
\boldsymbol{I}^{\text{RD,qs}},
\boldsymbol{U}^{\text{RD,qs}},\\
\boldsymbol{P}^{\text{RD}},
\boldsymbol{P}^{\text{RD,sg}},
\boldsymbol{S}^{\text{gs}},
\boldsymbol{S}^{\text{ls}},
\boldsymbol{S}^{+},
\boldsymbol{S}^{-},
\boldsymbol{W}^{\text{RD}}
     \end{array} \right\}                                   \mspace{-290mu}&                                        \notag\\
&\textit{Generator Constraints:}                \mspace{-290mu}& \notag       \\
&\textstyle{I_{it}^{\text{RD}}=I_{it}^{\star}+I_{it}^{\text{RD,qs}},}
                                                \mspace{-290mu}&\forall i\in\mathcal{I},t\in\mathcal{T};\label{DED:2}\\
&\textstyle{I_{it}^{\text{RD,qs}} \leq O_{it}^{\star},}
                                                \mspace{-290mu}&\forall i\in\mathcal{I},t\in\mathcal{T};\label{DED:3}\\
&\textstyle{U_{it}^{\text{RD,qs}} - D_{it}^{\text{RD,qs}}=I_{it}^{\text{RD,qs}}-I_{i,t-1}^{\text{RD,qs}},}
                                                \mspace{-290mu}&\forall i\in\mathcal{I},t\in\mathcal{T};\label{DED:4}\\
&\textstyle{P_{it}^{\text{RD}}=\sum\nolimits_{k \in \mathcal{K}}P_{itk}^{\text{RD,sg}},}
                                                \mspace{-290mu}&\forall i\in\mathcal{I},t\in\mathcal{T};\label{DED:5}\\
&\textstyle{0 \leq P_{itk}^{\text{RD,sg}} \leq \bar{P}_{ik}^{\text{sg}}I_{it}^{\text{RD}},}
                                                \mspace{-290mu}&\forall i\in\mathcal{I},t\in\mathcal{T},k\in\mathcal{K};
                                                                                                        \label{DED:6}\\
&\textstyle{P_{it}^{\text{RD}} \geq P^{\text{m}}_{i}I_{it}^{\text{RD}},}
                                                 \mspace{-290mu}&\forall i\in\mathcal{I},t\in\mathcal{T};\label{DED:7}\\
&\textstyle{P_{it}^{\text{RD}} \leq P^{\text{M}}_{i}I_{it}^{\star} + R_{it}^{\text{nr}\star}I_{it}^{\text{RD,qs}},}
                                                 \mspace{-290mu}&\forall i\in\mathcal{I},t\in\mathcal{T};\label{DED:8}\\
&\textstyle{P_{it}^{\text{RD}} - P_{it}^{\star} \geq -R_{it}^{\text{sr}\star} I_{it}^{\star} - {P}^{\text{M}}_{i}I_{it}^{\text{RD,qs}},}
                                                 \mspace{-290mu}&\forall i\in\mathcal{I},t\in\mathcal{T};\label{DED:9}\\
&\textstyle{P_{it}^{\text{RD}} - P_{it}^{\star} \leq  R_{it}^{\text{sr}\star} I_{it}^{\star} + {P}^{\text{M}}_{i}I_{it}^{\text{RD,qs}},}                                              \mspace{-290mu}&\forall i\in\mathcal{I},t\in\mathcal{T};\label{DED:10}\\
&\textstyle{P_{it}^{\text{RD}} - P_{i,t-1}^{\text{RD}} \leq P^{\text{M}}_{i}(1-I_{it}^{\text{RD}})}
                                                \mspace{-290mu}&\notag       \\
&\textstyle{\mspace{30mu} + R_{i}^{\uparrow} I_{i,t-1}^{\text{RD}} + R_{i}^{\text{su}} (I_{it}^{\text{RD}} - I_{i,t-1}^{\text{RD}})},                                    \mspace{-290mu}&\forall i\in\mathcal{I},t\in\mathcal{T};\label{DED:11}\\
&\textstyle{P_{i,t-1}^{\text{RD}} - P_{it}^{\text{RD}} \leq P^{\text{M}}_{i}(1-I_{i,t-1}^{\text{RD}})}
                                                \mspace{-290mu}&\notag       \\
&\textstyle{\mspace{30mu} + R^{\downarrow}_{i} I_{it}^{\text{RD}} + R^{\text{sd}}_{i}(I_{i,t-1}^{\text{RD}} - I_{it}^{\text{RD}}),}                                     \mspace{-290mu}&\forall i\in\mathcal{I},i\in\mathcal{I};\label{DED:12}\\
&\textstyle{0\leq W_{jt}^{\text{RD}}\leq\tilde{{W}}_{jt},}
                                                \mspace{-290mu}&\forall j\in\mathcal{J},t\in\mathcal{T};\label{DED:13}\\
&\textstyle{
D_{it}^{\text{RD,qs}},
I_{it}^{\text{RD}},
I_{it}^{\text{RD,qs}},
U_{it}^{\text{RD,qs}}\in\{0,1\}},
                                                \mspace{-290mu}&\forall i\in\mathcal{I},t\in\mathcal{T};\label{DED:14}\\
&\textit{System Constraints:}                   \mspace{-290mu}&\notag       \\
&\textstyle{    \sum\limits_{i \in \mathcal{I}} P_{it}^{\text{RD}}
               +\sum\limits_{j \in \mathcal{J}} W_{jt}^{\text{RD}}
               -                                S_{t}^{\text{gs}}
               =\sum\limits_{q \in \mathcal{Q}} \tilde{L}_{qt}
               - S_{t}^{\text{ls}},}            \mspace{-290mu}&\forall t\in\mathcal{T};                \label{DED:15}\\
&\textstyle{\mathcal{F}_{b}(\boldsymbol{P}_{t}^{\text{RD}}, \boldsymbol{W}_{t}^{\text{RD}},
\tilde{\boldsymbol{L}}_{t}) - S_{bt}^{-} \leq B_{b},}       \mspace{-290mu}&\forall b\in\mathcal{B},t\in\mathcal{T};\label{DED:16}\\
&\textstyle{\mathcal{F}_{b}(\boldsymbol{P}_{t}^{\text{RD}}, \boldsymbol{W}_{t}^{\text{RD}}, \tilde{\boldsymbol{L}}_{t}) + S_{bt}^{+} \geq -B_{b},}
                                                \mspace{-290mu}&\forall b\in\mathcal{B},t\in\mathcal{T};\label{DED:17}\\
&\textstyle{S_{t}^{\text{gs}},S_t^{\text{ls}}\geq 0,} \mspace{-290mu}&                \forall t\in\mathcal{T};\label{DED:18}\\
&\textstyle{S_{bt}^{+}, S_{bt}^{-} \geq 0,}     \mspace{-290mu}&\forall b\in\mathcal{B},t\in\mathcal{T};\label{DED:19}
\end{flalign}
\end{subequations}

Following the current industry practice, the DC power flow model is used in this paper to formulate network security constraints. Nevertheless, the presented bi-level MIP framework can incorporate other complex forms of network security constraints explored in the academic literature, such as linearized \cite{new} and SOCP-based \cite{acopf} AC power flow formulations.

\subsection{Relaxed Operation Models of CPO-NT-$\phi$ With Tunable Relaxing Degree}
\label{app_c}

CPO-NT-$\phi$ is trained by the bi-level methodology described in Section~\ref{sec03}, in which the binary requirement in the lower-level operation models is relaxed.

The relaxed UC and RD models are formulated as in \eqref{LPUC} and \eqref{LPRD}, in which the binary variables are relaxed as in \eqref{LPUC:3} and \eqref{LPRD:3}, while their relaxing degrees are further controlled by constraints \eqref{LPUC:4}-\eqref{LPUC:8} and \eqref{LPRD:4}-\eqref{LPRD:8} via a pre-determined parameter $\phi \in [\text{0}, \text{1}]$ together with four sets of new variables \cite{relax}, i.e., those hatted by $\prime$.
\begin{subequations}\label{LPUC}
\begin{flalign}
      &      \mspace{-30mu}\textit{$\phi$-Relaxed Unit Commitment Model:}        \mspace{-320mu} &\notag\\
\min  & \,\text{ The original UC objective function \eqref{DUC:1}}   \mspace{-320mu} &\label{LPUC:1}\\
s.\,t.&\,\eqref{DUC:2}-\eqref{DUC:15}, \eqref{DUC:17}-\eqref{DUC:19};\mspace{-320mu} &\label{LPUC:2}\\
      &\,0\leq D_{it}, I_{it}, O_{it}, U_{it}\leq 1,                 \mspace{-320mu} &\forall i\in\mathcal{I},t\in\mathcal{T};                                                                                \label{LPUC:3}\\
      &\,D_{it}=\textstyle{(1-\phi)\sum_{a=1}^{\lfloor \log_{2}\frac{1}{1-\phi} \rfloor + 1}}2^{a-1} \acute{D}_{ita},
                                                                     \mspace{-320mu}&\notag\\
      &                                                              \mspace{-320mu}&\forall i\in\mathcal{I},t\in\mathcal{T};                                                      \label{LPUC:4}\\
      &\,I_{it}=\textstyle{(1-\phi)\sum_{a=1}^{\lfloor \log_{2}\frac{1}{1-\phi} \rfloor + 1}}2^{a-1} \acute{I}_{ita},
                                                                     \mspace{-320mu} &\notag\\
      &                                                              \mspace{-320mu} &\forall i\in\mathcal{I},t\in\mathcal{T};                                                      \label{LPUC:5}\\
      &\,O_{it}=\textstyle{(1-\phi)\sum_{a=1}^{\lfloor \log_{2}\frac{1}{1-\phi} \rfloor + 1}}2^{a-1} \acute{O}_{ita},
                                                                     \mspace{-320mu}&\notag\\
      &                                                              \mspace{-320mu}&\forall i\in\mathcal{I},t\in\mathcal{T};                                                      \label{LPUC:6}\\
      &\,U_{it}=\textstyle{(1-\phi)\sum_{a=1}^{\lfloor \log_{2}\frac{1}{1-\phi}  \rfloor + 1}}2^{a-1} \acute{U}_{ita},
                                                                     \mspace{-320mu}& \notag\\
      &                                                              \mspace{-320mu}&\forall i\in\mathcal{I},t\in\mathcal{T};                                                      \label{LPUC:7}\\
      &\,\acute{D}_{ita}, \acute{I}_{ita}, \acute{O}_{ita}, \acute{U}_{ita}\in\{0, 1\},
                                                                     \mspace{-320mu}&\notag\\
      &                                                              \mspace{-320mu}&\forall i\in\mathcal{I},t\in\mathcal{T},a\in \{1,\cdots,\lfloor \log_{\text{2}}\frac{\text{1}}{\text{1}-\phi} \rfloor + \text{1}\};                                                                           \label{LPUC:8}
\end{flalign}
\end{subequations}
\begin{subequations}\label{LPRD}
\vspace{-5mm}
\begin{flalign}
       &\mspace{-30mu} \textit{$\phi$-Relaxed Re-Dispatch Model:}     \mspace{-320mu}&\notag\\
\min  & \,\text{ The original RD objective function \eqref{DED:1}}    \mspace{-320mu}&\label{LPRD:1}\\
s.\,t.&\,\eqref{DED:2}-\eqref{DED:13}, \eqref{DED:15}-\eqref{DED:19}; \mspace{-320mu}&\label{LPRD:2}\\
      &\,0 \leq D_{it}^{\text{RD,qs}}, I_{it}^{\text{RD}},I_{it}^{\text{RD,qs}},U_{it}^{\text{RD,qs}} \leq 1,
                                                                      \mspace{-320mu}&\forall  i \in \mathcal{I}, t \in \mathcal{T}; \label{LPRD:3}\\
&\, D_{it}^{\text{RD,qs}}=\textstyle{(1-\phi)\sum_{a=1}^{\lfloor \log_{2}\frac{1}{1-\phi}  \rfloor + 1}}2^{a-1} \acute{D}_{ita}^{\text{RD,qs}},                                         \mspace{-320mu}&\notag\\
      &                                                                 \mspace{-320mu}&\forall i \in \mathcal{I}, t \in \mathcal{T}; \label{LPRD:4}\\
      &\, I_{it}^{\text{RD}}=\textstyle{(1-\phi)\sum_{a=1}^{\lfloor \log_{2}\frac{1}{1-\phi}  \rfloor + 1}}2^{a-1} \acute{I}_{ita}^{\text{RD}},                                           \mspace{-320mu}&\notag\\
      &                                                                \mspace{-320mu}&\forall  i \in \mathcal{I}, t \in \mathcal{T}; \label{LPRD:5}\\
      &\, I_{it}^{\text{RD,qs}}=\textstyle{(1-\phi)\sum_{a=1}^{\lfloor \log_{2}\frac{1}{1-\phi} \rfloor + 1}}2^{a-1} \acute{I}_{ita}^{\text{RD,qs}},                                        \mspace{-320mu}&\notag\\
      &                                                                \mspace{-320mu}&\forall i \in \mathcal{I}, t \in \mathcal{T}; \label{LPRD:6}\\
      &\, U_{it}^{\text{RD,qs}}=\textstyle{(1-\phi)\sum_{a=1}^{\lfloor \log_{2}\frac{1}{1-\phi}  \rfloor + 1}}2^{a-1} \acute{U}_{ita}^{\text{RD,qs}},                                         \mspace{-320mu}&\notag\\
      &                                                                 \mspace{-320mu}&\forall i \in \mathcal{I}, t \in \mathcal{T}; \label{LPRD:7}\\
      &\,\acute{D}_{ita}^{\text{RD,qs}}, \acute{I}_{ita}^{\text{RD}},\acute{I}_{ita}^{\text{RD,qs}},\acute{U}_{ita}^{\text{RD,qs}} \in \{0, 1\};                                           \mspace{-320mu}&\notag\\
      &                                                                 \mspace{-320mu}&\forall i \in \mathcal{I}, t \in \mathcal{T}, a \in \{1,\cdots,\lfloor \log_{\text{2}}\frac{\text{1}}{\text{1}-\phi} \rfloor + \text{1}\}; \label{LPRD:8}
\end{flalign}
\end{subequations}

Here, $\lfloor \cdot \rfloor$ is the floor function. Models \eqref{LPUC} and \eqref{LPRD} highlight that their relaxing degree can be tuned via parameter $\phi$. Specifically, taking the status variable $I_{it}$ as an example, if $\phi=\text{0}$ (0\% relaxation), then $I_{it} \in \{\text{0}, \text{1}\}$\textemdash \eqref{LPUC} and \eqref{LPRD} are equivalent to the original MIP models \eqref{DUC} and \eqref{DED}; if $\phi=\text{0.5}$ (50\% relaxation), then $I_{it} \in \{\text{0}, \text{0.5}, \text{1}\}$; if $\phi=\text{0.9}$ (90\% relaxation), then $I_{it} \in \{\text{0}, \text{0.1}, \text{0.2}, \cdots, \text{0.9}, \text{1}\}$. Note that when $\phi=\text{1}$ (100\% relaxation), \eqref{LPUC:4}-\eqref{LPUC:8} and \eqref{LPRD:4}-\eqref{LPRD:8} will be naturally removed, transforming \eqref{LPUC} and \eqref{LPRD} into LP.

Regarding the training of CPO-NT-$\phi$, when $\phi \in [\text{0}, \text{1})$, \eqref{LPUC} and \eqref{LPRD} remain MIP models, and the training method described in Section~\ref{Solving} is still valid; when $\phi=\text{1}$, \eqref{LPUC} and \eqref{LPRD} become LP models, thus, the training follows the method presented in \cite{vladimir}\textemdash replacing the lower-level problems with their KKT conditions, then solving the ERM as an MIP.
\begin{figure}[H]
	\centering
		\includegraphics[width=0.92\columnwidth]{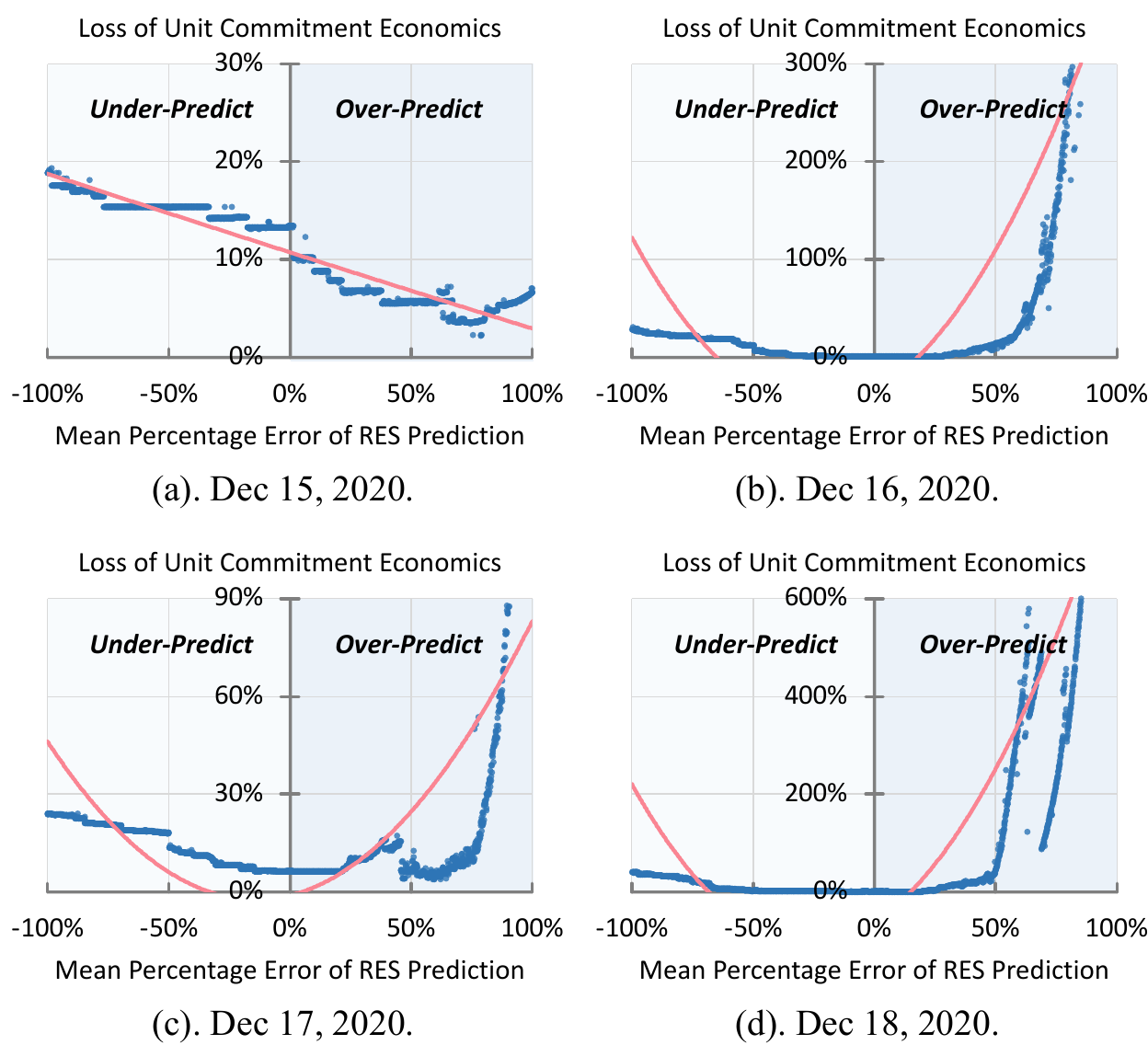}
		\vspace{-2mm}
	\caption{Approximating results for NNPO-NT.}\label{fig12}
\end{figure}

\begin{table}[H]
\vspace{-4mm}
	\caption{Coefficients of Approximating Results}\label{tab12}
	\centering
	\footnotesize
\begin{tabularx}{\columnwidth}{@{\extracolsep{\fill}} rrrr}
\toprule
             &\multicolumn{1}{c}{$a$}  &\multicolumn{1}{c}{$b$}          &\multicolumn{1}{c}{$c$}\\
\midrule
Dec 15, 2020 & 0.0012                  &-0.0791                          & 0.1074\\
Dec 16, 2020 & 2.9703                  & 1.4027                          &-0.3433\\
Dec 17, 2020 & 0.6565                  & 0.1846                          &-0.0097\\
Dec 18, 2020 & 6.0497                  & 3.2362                          &-0.6100\\
\bottomrule
\end{tabularx}
\end{table}

\begin{table}[H]
\vspace{-4mm}
	\caption{Hyper-Parameters of NNPO-NT}\label{tab13}
	\centering
	\footnotesize
\begin{tabularx}{\columnwidth}{@{\extracolsep{\fill}} ccc}
\toprule
                         &NNPO-7                   &NNPO-349\\
\midrule
Learning Rate            &0.001                    &0.0001   \\
Hidden Layer Structure   &[$|\mathcal{T}|\times|\mathcal{J}|$]           &[$|\mathcal{T}|\times|\mathcal{J}|$]\\
Activation Function      &linear                   &linear \\
Number of Epochs         &50                       &400\\
Batch Size               &1                        &2\\
\bottomrule
\end{tabularx}
\end{table}

\subsection{Implementation Details of NNPO-NT}\label{app_d}
NNPO-NT is a variant of the cost-oriented prediction method \cite{yiwang1}, including the following major steps:

\textit{Step 1:} For day $D$, 2,000 scenarios are generated to simulate a set of imperfect RES predictions. Then, the UC-RD process is executed for each scenario to calculate its UC economics loss \eqref{EconomicsLoss}. This procedure yields 2,000 data points describing the relationship between the mean percentage error (MPE) and UC economics loss, as depicted by the blue dots in Fig.~\ref{fig12};

\textit{Step 2:} Approximate the blue dots using quadratic functions \eqref{appro}, as depicted by the red lines in Fig.~\ref{fig12}. These quadratic functions approximate the relationship between MPE and UC economics loss \eqref{EconomicsLoss}. Coefficients of the approximate quadratic functions are listed in Table~\ref{tab12}. Note that the loss function \eqref{appro} can be linear, quadratic, cubic, or quartic. The quadratic function is selected for two reasons: \textit{i)} the insights gained from \cite{yiwang1} and \textit{ii)} its outperformance over other options.
\begin{equation}\label{appro}
\text{Loss Caused by MPE} = a\times\text{MPE}^{2} + b\times\text{MPE}^{1} + c
\end{equation}

\textit{Step 3:} Using the approximate functions \eqref{appro} as loss functions, the NN-based tailors are trained via TensorFlow 2.14. The setting of the hyper-parameters is listed in Table~\ref{tab13}.

\bibliographystyle{IEEEtran}
\bibliography{ref_CPO}

\begin{IEEEbiographynophoto}{Xianbang Chen} (Student Member, IEEE) received the B.S. and M.S. degrees in electrical engineering from Sichuan University, Chengdu, China, in 2017 and 2020, respectively. He is currently working toward the Ph.D. degree with the Stevens Institute of Technology, Hoboken, NJ, USA. His research interest focuses on applying prescription methods on unit commitment. \end{IEEEbiographynophoto}

\begin{IEEEbiographynophoto}{Yikui Liu} (Member, IEEE) received the B.S. degree in electrical engineering and automation from the Nanjing Institute of Technology, China, in 2012, the M.S. degree in power system and automation from Sichuan University, China, in 2015, and the Ph.D. degree in electrical and computer engineering from the Stevens Institute of Technology, Hoboken, NJ, USA, in 2020. He worked in Siemens, USA, as an Energy Market Engineer, during 2020–2021. He worked as a Postdoctoral Researcher with the Stevens Institute of Technology, Hoboken, NJ, USA, during 2021–2023. He is currently an associate researcher with Sichuan University, Chengdu China. His research interests include the power market and OPF in distribution systems.\end{IEEEbiographynophoto}

\begin{IEEEbiographynophoto}{Lei Wu} (Fellow, IEEE) received the B.S. degree in electrical engineering and the M.S. degree in systems engineering from Xi’an Jiaotong University, Xi’an, China, in 2001 and 2004, respectively, and the Ph.D. degree in electrical engineering from Illinois Institute of Technology (IIT), Chicago, IL, USA, in 2008. From 2008 to 2010, he was a Senior Research Associate with the Robert W. Galvin Center for Electricity Innovation, IIT. He was a summer Visiting Faculty at NYISO in 2012. He was a Professor with the Electrical and Computer Engineering Department, Clarkson University, Potsdam, NY, USA, till 2018. Currently, he is Anson Wood Burchard Chair Professor with the Department of Electrical and Computer Engineering, Stevens Institute of Technology, Hoboken, NJ, USA. His research interests include power systems operation and planning, energy economics, and community resilience microgrid.\end{IEEEbiographynophoto}

\end{document}